\documentclass[11pt,english]{amsart}
\usepackage[T1]{fontenc}
\usepackage[latin1]{inputenc}
\usepackage{a4wide}
\usepackage{amssymb}

 \theoremstyle{plain}    
 \newtheorem{thm}{Theorem}[section]
 \numberwithin{equation}{section} 
 \numberwithin{figure}{section} 
 \theoremstyle{definition}
 \newtheorem{defn}[thm]{Definition}
 \theoremstyle{remark}    
 \newtheorem{notation}[thm]{Notation} 
 \theoremstyle{remark}
 \newtheorem{rem}[thm]{Remark}
 \theoremstyle{definition}
  \newtheorem{example}[thm]{Example}
 \theoremstyle{plain}    
 \newtheorem{prop}[thm]{Proposition} 
 \theoremstyle{plain}    
 \newtheorem{lem}[thm]{Lemma} 
 \theoremstyle{plain}    
 \newtheorem{cor}[thm]{Corollary} 

\usepackage{latexsym,amsfonts, xypic,amssymb}
\usepackage{pstricks,pst-node}
\xyoption{all}
\newcommand{\bc}{\mathbb{C}}
\newcommand{\bcp}{\mathbb{C}_+}

\newcommand{\ba}{\mathbb{A}}

\newcommand{\gaz}{\mathbb{G}_{a,Z}}
\newcommand{\der}{\partial}
\theoremstyle{remark}

\theoremstyle{definition}
\newtheorem{enavant}[thm]{}

\usepackage{babel}
\begin{document}

\title{Embeddings of generalized Danielewski surfaces in affine spaces}

\maketitle
\begin{center}\author{{\bf A. DUBOULOZ} \\ \vspace{0.8cm}   \address Institut Fourier, Laboratoire de Mathématiques\\ UMR5582 (UJF-CNRS)\\ BP 74, 38402 ST MARTIN D'HERES CEDEX\\ FRANCE \\ \email adrien.dubouloz@ujf-grenoble.fr \\ \vspace{1cm}} \end{center}

\begin{abstract}
We construct explicit embeddings of generalized Danielewski surfaces
\cite{DubGDS} in affine spaces. The equations defining these embeddings
are obtain from the $2\times2$ minors of a matrix attached to a labelled-rooted
tree $\mathfrak{g}$. The corresponding surfaces $V_{\mathfrak{g}}$
come equipped with canonical $\bcp$-actions which appear as the restrictions
of certain $\bcp$-actions on the ambient affine space. We characterize
those surfaces $V_{\mathfrak{g}}$ with a trivial Makar-Limanov invariant,
and we complete the study of log $\mathbb{Q}$-homology planes with
a trivial Makar-Limanov invariant initiated by Miyanishi and Masuda
\cite{MiMa03}. 
\end{abstract}

\section*{Introduction}

In \cite{DubGDS}, the author introduced the notion of generalized
Danielewski surface ($GDS$ for short). This is a normal affine surface
$V$ which admits a faithfully flat morphism $q:V\rightarrow Z=\ba_{\bc}^{1}$
with reduced fibers and generic fiber $\ba_{K\left(Z\right)}^{1}$,
such that all but possibly one closed fiber $q^{-1}\left(z_{0}\right)$
are integral. For instance, a nonsingular ordinary Danielewski surface
$V_{P,n}\subset\bc\left[x,y,z\right]$ with equation $x^{n}z-P\left(y\right)=0$,
where $P$ is a nonconstant polynomial with simple roots, is a $GDS$
for the projection $pr_{x}:V_{P,n}\rightarrow Z=\textrm{Spec}\left(\bc\left[x\right]\right)$.
Generalized Danielewski surfaces appear naturally as locally trivial
fiber bundles $\rho:V\rightarrow X$ over an affine line with a multiple
origin (see e.g. \cite{Fie94}). In \cite{DubGDS}, the author established
that every such bundle is obtain from an invertible sheaf $\mathcal{L}$
on $X$ and a \v{C}ech $1$-cocycle $g$ with values in the dual
$\mathcal{L}^{\vee}$ of $\mathcal{L}$. In turn, this invertible
sheaf $\mathcal{L}$ and the cocycle $g$ are encoded in a combinatorial
data consisting of a rooted tree with weighted edges, which we call
a \emph{weighted rooted tree} (see \cite[Theorem 4.1]{DubGDS} and
\ref{GDS2WTree} below). Here we use rooted trees in a different way
to construct embeddings of $GDS$'s into affine spaces. More precisely,
starting with certain rooted trees with weights on their nodes, which
we call \emph{labelled rooted trees}, we construct explicit ideals
of certain polynomial rings. In turn, these ideals define affine surfaces
which are naturally $GDS$'s over the affine line $\ba_{\bc}^{1}$.

\psset{unit=1.2cm}

\begin{pspicture}(-2,-2)(6,2)

\psline(0,0)(1,0)(2,-0.5)(3,-0.5) \psline(1,0)(2,0.5)

\rput(0,0){\textbullet}\rput(1,0){\textbullet}\rput(2,-0.5){\textbullet}\rput(2,0.5){\textbullet}\rput(3,-0.5){\textbullet}

\rput(-0.1,-0.4){{\small $e_0$}}

\rput(0.5,0.2){{\small $0$}}

\rput(1.3,0.4){{\small $1$}}\rput(1.2,-0.4){{\small $-1$}}

\rput(2.5,-0.3){{\small $0$}}

\rput(1,-1.5){\textrm{ A weighted tree} $\gamma$ \textrm{rooted in} $e_0$.} 

\psline(5,0)(7,0)(8,0.5)\psline(5,0)(6,0.5) \psline(5,0)(6,-0.5)\psline(7,0)(8,-0.5)

\rput(5,0){\textbullet}\rput(6,0){\textbullet}\rput(6,0.5){\textbullet}\rput(6,-0.5){\textbullet}\rput(7,0){\textbullet}\rput(8,0.5){\textbullet}\rput(8,-0.5){\textbullet}

\rput(4.9,-0.4){{\small $e_0$}}

\rput(6,0.2){{\small $0$}}

\rput(6,0.8){{\small $1$}}

\rput(6,-0.8){{\small $2$}}

\rput(7,0.2){{\small $3$}}

\rput(8,0.8){{\small $0$}}

\rput(8,-0.8){{\small $0$}}

\rput(6.5,-1.5){\textrm{ A labelled tree} $\mathfrak{g}$ \textrm{rooted in} $e_0$.} 

\end{pspicture}

The paper is divided as follows. In section $1$ we recall basic facts
on rooted trees and we give a procedure to construct an affine scheme
$V_{\mathfrak{g}}$ from the data consisting of a labelled rooted
tree $\mathfrak{g}=\left(\Gamma,lb\right)$. These schemes $V_{\mathfrak{g}}=\textrm{Spec}\left(B_{\mathfrak{g}}\right)$
are naturally schemes over an affine line $Z=\textrm{Spec}\left(\bc\left[h\right]\right)$,
and they come embedded into affine spaces $\ba_{Z}^{d\left(\Gamma\right)+1}=\textrm{Spec}\left(\bc\left[h\right]\left[\Gamma\right]\right)$
depending on the tree $\Gamma$. In section $2$, we prove the following
result (theorem \ref{EmbeddedGDSTheorem})

\begin{thm}
For every fine-labelled rooted tree $\mathfrak{g}$ (see definition
\ref{LabelDef} below), the scheme $V_{\mathfrak{g}}=\textrm{Spec}\left(B_{\mathfrak{g}}\right)$
is a $GDS$ over $Z=\textrm{Spec}\left(\bc\left[h\right]\right)$. 
\end{thm}
\noindent For instance, the Bandman and Makar-Limanov surface \cite{BML01}
$V\subset\bc\left[x,y,z,u\right]$ with equation \[
xz=y\left(y^{2}-1\right),\quad yu=z\left(z^{2}-1\right),\quad xu=\left(y^{2}-1\right)\left(z^{2}-1\right)\]
 is a $GDS$ over $Z=\textrm{Spec}\left(\bc\left[x\right]\right)$
for the morphism $pr_{x}:V\rightarrow Z$. This data corresponds to
the fine-labelled rooted tree

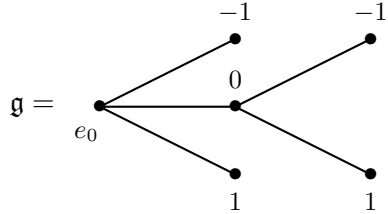
\psset{unit=1.8cm}

\begin{pspicture}(-3,-1)(3,1)

\rput(-0.5,0){$\mathfrak{g}=$}

\psline(0,0)(1,0.5)\psline(0,0)(1,0)\psline(0,0)(1,-0.5)

\psline(1,0)(2,0.5)\psline(1,0)(2,-0.5)

\rput(0,0){\textbullet}

\rput(1,0){\textbullet}\rput(1,0.5){\textbullet}\rput(1,-0.5){\textbullet}

\rput(2,0.5){\textbullet}\rput(2,-0.5){\textbullet}

\rput(-0.1,-0.2){{\small $e_0$}}

\rput(1,0.7){{\small $-1$}}

\rput(1,-0.7){{\small $1$}}

\rput(1,0.2){{\small $0$}}

\rput(2,0.7){{\small $-1$}}

\rput(2,-0.7){{\small $1$}}

\end{pspicture}

\noindent In section $2$, we recall \cite{DubGDS} how a weighted
rooted tree $\gamma$ defines a $GDS$ $q_{\gamma}:W_{\gamma}\rightarrow Z$.
Then, starting with such a weighted rooted tree $\gamma=\left(\Gamma,w\right)$,
we construct a fine-labelled rooted tree $\mathfrak{g}_{w}=\left(\Gamma,lb_{w}\right)$
with the same underlying tree $\Gamma$ and a closed embedding $i_{\gamma}:W_{\gamma}\hookrightarrow\ba_{Z}^{d\left(\Gamma\right)+1}$
inducing an isomorphism $W_{\gamma}\simeq V_{\mathfrak{g}_{w}}$.
This leads to the following result (theorem \ref{Flrt2GDSClassifTheorem})

\begin{thm}
Every $GDS$ $q:V\rightarrow Z$ is $Z$-isomorphic to a $GDS$ $q_{\mathfrak{g}}:V_{\mathfrak{g}}\rightarrow Z$
for an appropriate fine-labelled rooted tree $\mathfrak{g}$.
\end{thm}
\noindent In section $3$ we study $\gaz$-actions on the $GDS$'s
$q_{\mathfrak{g}}:V_{\mathfrak{g}}=\textrm{Spec}\left(B_{\mathfrak{g}}\right)\rightarrow Z$.
We construct explicit locally nilpotent derivations of the algebras
$B_{\mathfrak{g}}$. In \cite{ML90,ML01}, Makar-Limanov proved that
every $\bcp$-action on an ordinary Danielewski surface $V_{P,n}\subset\textrm{Spec}\left(\bc\left[x,y,z\right]\right)$
is the restriction of a $\bcp$-action on $\ba_{\bc}^{3}=\textrm{Spec}\left(\bc\left[x,y,z\right]\right)$.
We prove the following general result (theorem \ref{ActionExtensionProp}). 

\begin{thm}
Every $\gaz$-action on a $GDS$ $q:V_{\mathfrak{g}}\rightarrow Z$
defined by a fine-labelled rooted tree $\mathfrak{g}=\left(\Gamma,lb\right)$
is induced by a $\gaz$-action on the ambient space $\ba_{Z}^{d\left(\Gamma\right)+1}=\textrm{Spec}\left(\bc\left[h\right]\left[\Gamma\right]\right)$.
\end{thm}
\noindent In section $4$, we consider $GDS$'s $q_{\mathfrak{g}}:V_{\mathfrak{g}}\rightarrow Z$
with a trivial Makar-Limanov invariant. As a consequence of \cite[Theorem 7.2]{DubGDS},
we obtain the following criterion (theorem \ref{Comb2GDSThm}).

\begin{thm}
A $GDS$ $q_{\mathfrak{g}}:V_{\mathfrak{g}}\rightarrow Z$ has a trivial
Makar-Limanov invariant if and only if $\mathfrak{g}$ is a fine-labelled
comb (see \ref{CombNotation} below).
\end{thm}
\noindent This leads to the following description (see \ref{NewEquations}
below). Let $P_{1},\ldots,P_{n}\in\bc\left[T\right]$ be a collection
of nonconstant polynomials with simple roots, one of these root, say
$\lambda_{i,1}$, $1\leq i\leq n$, being distinguished. We let \[
P_{i}\left(T\right)=\prod_{j=1}^{r_{i}}\left(T-\lambda_{i,j}\right)=\left(T-\lambda_{i,1}\right)\tilde{P}_{i}\left(T\right)\in\bc\left[T\right],\quad1\leq i\leq n.\]
 Then every $GDS$ with a trivial Makar-Limanov is isomorphic, for
an appropriate choice of the polynomials $P_{1},\ldots,P_{n}$, to
a surface \[
V_{P_{1},\ldots,P_{n}}\subset\ba_{\bc}^{n+2}=\textrm{Spec}\left(\bc\left[x\right]\left[y_{1},\ldots,y_{n+1}\right]\right)\]
 with equations \[
\left\{ \begin{array}{lcll}
xy_{i+1} & = & \left({\displaystyle \prod_{k=1}^{i-1}}\tilde{P}_{k}\left(y_{k}\right)\right)P_{i}\left(y_{i}\right) & 1\leq i\leq n\\
\left(y_{j-1}-\lambda_{j-1,1}\right)y_{i+1} & = & y_{j}\left({\displaystyle \prod_{k=j}^{i-1}}\tilde{P}_{k}\left(y_{k}\right)\right)P_{i}\left(y_{i}\right) & 2\leq j\leq i\leq n\end{array}\right.,\]
 where, by convention, ${\displaystyle {\displaystyle \prod_{k=j}^{i-1}}}\tilde{P}_{k}\left(y_{k}\right)=1$
if $i=j$. From this description, we recover the following characterization
of nonsingular ordinary Danielewski surfaces, due to Bandman and Makar-Limanov
\cite{BML01} (\ref{MLGDSThm}).

\begin{thm}
For a $GDS$ $q:V\rightarrow Z$ with a trivial Makar-Limanov invariant,
the following are equivalent.

1) $V$ admits a free $\gaz$-action.

2) The canonical sheaf $\omega_{V}$ of $V$ is trivial.

3) $V$ is isomorphic to an ordinary Danielewski surface $V_{P,1}\subset\ba_{\bc}^{3}=\textrm{Spec}\left(\bc\left[x,y,z\right]\right)$
with the equation $xz-P\left(y\right)=0$ for a certain nonconstant
polynomial $P$ with simple roots.
\end{thm}
\noindent Following a remark of the author, Miyanishi and Masuda
\cite{MiMa03} proved that every nonsingular $\mathbb{Q}$-homology
plane with a trivial Makar-Limanov invariant is a cyclic quotient
of an ordinary Danielewski surface $V\subset\textrm{Spec}\left(\bc\left[x,t,z\right]\right)$
with equation $xz=t^{m}-1$ by a $\mathbb{Z}_{m}$-action $\left(x,t,z\right)\mapsto\left(\varepsilon x,\varepsilon^{q}y,\varepsilon^{-1}z\right)$,
where $\varepsilon$ is a primitive $m$-the root of unity and $\gcd\left(q,m\right)=1$.
More generally, the author proved in \cite[theorem 7.7]{DubGDS} that
every normal affine surface $S$ with a trivial Makar-Limanov invariant
is a cyclic quotient of a $GDS$ $V$. In case that $S$ is a log
$\mathbb{Q}$-homology plane, this $GDS$ can be determined explicitly.
This leads to the following result (theorem \ref{LogQHomThm}), see
also \cite{DaRus02}.

\begin{thm}
Every log $\mathbb{Q}$-homology plane $S\not\simeq\ba_{\bc}^{2}$
with a trivial Makar-Limanov invariant is isomorphic to a quotient
of an ordinary Danielewski surface $V\subset\textrm{Spec}\left(\bc\left[x,t,z\right]\right)$
with equation $xz=t^{n}-1$ by a $\mathbb{Z}_{m}$-action $\left(x,t,z\right)\mapsto\left(\varepsilon x,\varepsilon^{q}t,\varepsilon^{-1}z\right)$,
where $\varepsilon$ is a primitive $m$-th root of unity, $n$ divides
$m$ and $\gcd\left(q,m/n\right)=1$. 
\end{thm}

\section{Affine schemes defined from labelled rooted trees}

In this section we explain how to construct affine schemes from labelled
rooted trees.

\subsection{Basics on rooted trees}

\indent

\noindent Let $G=\left(N,\leq\right)$ be a nonempty, finite, partially
ordered set (a \emph{poset}, for short). The elements of $N$ are
sometimes called the \emph{nodes} of $G$. A totally ordered subset
$N'\subset N$ is called a \emph{chain of length} $l\left(N'\right)=\textrm{Card}\left(N'\right)-1$.
\emph{A} chain which is maximal for the inclusion is called a \emph{maximal
chain.} For every $e\in N$, we let \[
\left(\uparrow e\right)_{G}=\left\{ e'\in N,e'\geq e\right\} \quad\textrm{and}\quad\left(\downarrow e\right)_{G}=\left\{ e'\in N,e'\leq e\right\} \]
 For every $e<e'$, $e,e'\in N$, we let $\left[e',e\right]_{G}:=\left(\uparrow e'\right)_{G}\cap\left(\downarrow e\right)_{G}$.
A pair $e<e'$ such that $\left[e',e\right]_{G}=\left\{ e'<e\right\} $
is called an edge of $G$, and we denote the set of all edges in $G$
by $E\left(G\right)$. 

\begin{defn}
A \emph{rooted tree $\Gamma$} is poset $\Gamma=\left(N\left(\Gamma\right),\leq\right)$
with a unique minimal element $e_{0}$ called the \emph{root}, and
\emph{}such that $\left(\downarrow e\right)_{\Gamma}$ is a chain
for every $e\in N\left(\Gamma\right)$. 
\end{defn}
\begin{enavant} The maximal elements of a rooted tree $\Gamma=\left(N\left(\Gamma\right),\leq\right)$
are called the \emph{leaves} of $\Gamma$. We denote the set of those
elements by $\textrm{Leaves}\left(\Gamma\right)$. An element of $N\left(\Gamma\right)\setminus\textrm{Leaves}\left(\Gamma\right)$
is called a \emph{parent}, and we denote the set of those nodes by
$\mathbf{P}\left(\Gamma\right)$. Given $e\in N\left(\Gamma\right)\setminus\left\{ e_{0}\right\} $,
an element of the chain $\textrm{Anc}\left(e\right)=\left(\downarrow e\right)\setminus\left\{ e\right\} $
is called an \emph{ancestor} of $e$. The \emph{parent of $e$} is
the maximal element $\textrm{Par}\left(e\right)$ of $\textrm{Anc}\left(e\right)$.
More generally, the $n$-th \emph{ancestor} of $e$ is defined recursively
by $\textrm{Par}^{n}\left(e\right)=\textrm{Par}\left(\textrm{Par}^{n-1}\left(e\right)\right)\in\textrm{Anc}\left(e\right)$.
Given two different nodes $g,g'\in N\left(\Gamma\right)$, the \emph{first
common ancestor} of $g$ \emph{and} $g'$ is the maximal element $\textrm{Anc}\left(g,g'\right)$
of the chain $\textrm{Anc}\left(g\right)\cap\textrm{Anc}\left(g\right)$.
Given $e\in\mathbf{P}\left(\Gamma\right)$, the minimal elements of
$\left(\uparrow e\right)_{\Gamma}\setminus\left\{ e\right\} $ are
called the \emph{children} of $e$, and we denote the set of those
nodes by $\textrm{Child}_{\Gamma}\left(e\right)$. The \emph{degree}
$\deg_{\Gamma}\left(e\right)$ of a node $e$ is the number of its
children. Given $e\in N\left(\Gamma\right)$, the \emph{maximal rooted
subtree of} $\Gamma$ \emph{rooted in} $e$ is the tree $\Gamma\left(e\right)=\left(\left(\uparrow e\right)_{\Gamma},\leq\right)$.
A node $e\in N\left(\Gamma\right)$ such that $l\left(\left(\downarrow e\right)\right)=n$
is said to be at \emph{level} $n$, and we denote the set of those
nodes by $N_{n}\left(\Gamma\right)$. The maximal chains of a rooted
tree $\Gamma$ are the chains \begin{equation}
\left(\downarrow e\right)_{\Gamma}=\left\{ e_{0}<e_{1}<\cdots<e_{n_{f}-1}<e_{n_{f}}=e\right\} ,\quad e\in\textrm{Leaves}\left(\Gamma\right).\label{MaximalChainNotation}\end{equation}
The \emph{height} $h\left(\Gamma\right)$ of $\Gamma$ is the maximum
of the lengths $l\left(\left(\downarrow e\right)_{\Gamma}\right)=n_{e}$,
$e\in\textrm{Leaves}\left(\Gamma\right)$. 

\end{enavant}

\noindent In \ref{LabelDef} and \ref{WeightDef} below, we introduce
two different notions of weights on a rooted tree $\Gamma$. 

\begin{defn}
\label{LabelDef} A \emph{labelling} (with values in $\bc$) on a
tree $\Gamma$ rooted in $e_{0}$ is a function\[
lb:N\left(\Gamma\right)\setminus\left\{ e_{0}\right\} \rightarrow\bc\]
 A labelling $lb$ is called \emph{fine} if $lb\left(f\right)\neq lb\left(f'\right)$
whenever $f$ and $f'$ share the same parent $e\in\mathbf{P}\left(\Gamma\right)$.
A rooted tree $\Gamma$ equipped with a labelling $lb$ will be referred
to as a \emph{labelled rooted tree} and will be denoted by $\mathfrak{g}=\left(\Gamma,lb\right)$.
A rooted tree $\Gamma$ equipped with a fine labelling $lb$ will
be referred to as a \emph{fine-labelled rooted tree.} For every $e\in\mathbf{P}\left(\Gamma\right)$,
we let $\mathfrak{g}\left(e\right)=\left(\Gamma\left(e\right),lb\right)$
be the maximal subtree of $\Gamma$ rooted in $e$, equipped with
the restriction of the labelling $lb$ to $N\left(\Gamma\left(e\right)\right)\setminus\left\{ e\right\} $. 
\end{defn}
\begin{defn}\label{WeightDef} A \emph{weight function} (with values
in $\bc$) on a rooted tree $\Gamma$ is a function\[
w:E\left(\Gamma\right)\rightarrow\bc,\]
 i.e. a function which assign a complex number $w\left(\left[e,f\right]_{\Gamma}\right)$
to every edge $\left[e,f\right]_{\Gamma}$ of $\Gamma$, such that
$w\left(\left[e,f\right]\right)\neq w\left(\left[e,f'\right]\right)$
whenever $f$ and $f'$ share the same parent $e\in\mathbf{P}\left(\Gamma\right)$.
A rooted tree $\Gamma$ equipped with a weight function $w$ will
be referred to as a \emph{weighted rooted tree} and will be denoted
by $\gamma=\left(\Gamma,w\right)$\emph{.} For every $e\in\mathbf{P}\left(\Gamma\right)$,
we let $\gamma\left(e\right)=\left(\Gamma\left(e\right),w\right)$
be the maximal subtree of $\Gamma$ rooted in $e$, equipped with
the restriction of the weight function $w$ to $E\left(\Gamma\left(e\right)\right)$. 

\end{defn}

\subsection{Affine scheme associated to a labelled rooted tree}

\indent 

\noindent In this subsection we give a general procedure for constructing
an affine scheme $V_{\mathfrak{g}}$ from the data consisting of a
labelled rooted tree $\mathfrak{g}=\left(\Gamma,lb\right)$. To fix
the notation, we let $A=\bc\left[h\right]$ be a polynomial ring in
one variable, $z_{0}=\left(h\right)\in Z=\textrm{Spec}\left(A\right)$
and $Z_{*}=Z\setminus\left\{ z_{0}\right\} \simeq\textrm{Spec}\left(A_{h}\right)$.
We let $pr_{Z}:\ba_{Z}^{1}=\textrm{Spec}\left(A\left[X_{0}\right]\right)\rightarrow Z$
be the trivial line bundle over $Z$. 

\begin{defn}
\label{BaseDefinition} Given a rooted tree $\Gamma$, we let $\mathbf{S}\left(M_{\Gamma}\right)$
be the symmetric algebra of the free $A$-module $M_{\Gamma}$ with
basis ${\displaystyle \left(X_{e}\right)_{e\in\mathbf{P}\left(\Gamma\right)}}$.
This is a polynomial ring over $A$ in \begin{equation}
d\left(\Gamma\right):=\sum_{i=1}^{h\left(\Gamma\right)}\textrm{Card}\left(N_{i-1}\left(\Gamma\right)\setminus\textrm{Leaves}\left(\Gamma\right)\right)\label{AmbientSpaceDimension}\end{equation}
 variables. Then we let $A\left[\Gamma\right]=A\left[X_{0}\right]\otimes_{A}\mathbf{S}\left(M_{\Gamma}\right)$. 
\end{defn}
\begin{notation}
If $e'\in\mathbf{P}\left(\Gamma\right)$ is the parent of a given
$e\in\mathbf{P}\left(\Gamma\right)\setminus\left\{ e_{0}\right\} $
then we will sometimes denote $X_{e'}\in A\left[\Gamma\right]$ as
$X_{\textrm{Par}\left(e\right)}$. We also extend this relationship
between the variables $X_{e}$, $e\in\mathbf{P}\left(\Gamma\right)$,
by letting $X_{\textrm{Par}\left(e_{0}\right)}=X_{0}\in A\left[\Gamma\right].$ 
\end{notation}
\noindent For every node $e\in\mathbf{P}\left(\Gamma\right)$ of
a given labelled rooted tree $\mathfrak{g}=\left(\Gamma,lb\right)$,
we introduce, in \ref{SPolynomialDef}-\ref{QPolyDefinition} below,
three polynomials $S_{e}\left(\mathfrak{g}\right),R_{e}\left(\mathfrak{g}\right),Q_{e}\left(\mathfrak{g}\right)\in A\left[\Gamma\right]$,
defined recursively through the labelling $lb$. 

\begin{defn}
\label{SPolynomialDef} For every $e\in\mathbf{P}\left(\Gamma\right)$
and every subset $J\subset\textrm{Child}\left(e\right)$ we let \[
S_{e}^{J}\left(\mathfrak{g}\right)={\displaystyle \prod_{e''\in\left(\textrm{Child}\left(e\right)\setminus J\right)}}\left(X_{\textrm{Par}\left(e\right)}-lb\left(e''\right)\right)\in\bc\left[X_{\textrm{Par}\left(e\right)}\right]\subset A\left[\Gamma\right].\]

\noindent We call $S_{e}\left(\mathfrak{g}\right):={\displaystyle S_{e}^{\emptyset}\left(\mathfrak{g}\right)}$
\emph{the sibling polynomial} of $e$. This is a polynomial of degree
$\deg_{\Gamma}\left(e\right)$. The roots of $S_{e}\left(\mathfrak{g}\right)$
are simple if and only if $lb:N\left(\Gamma\right)\setminus\left\{ e_{0}\right\} \rightarrow\bc$
restricts to a fine labelling on the subtree $\left(\left\{ e\right\} \cup\textrm{Child}\left(e\right),\leq_{\Gamma}\right)$
of $\Gamma$. 
\end{defn}
\begin{defn}\label{RPolyDefinition} The \emph{root polynomial} of
$e\in N\left(\Gamma\right)$ is the polynomial $R_{e}\left(\mathfrak{g}\right)$
defined recursively by $R_{e_{0}}\left(\mathfrak{g}\right)=1$ and
\[
R_{e}\left(\mathfrak{g}\right)=S_{\textrm{Par}\left(e\right)}^{\left\{ e\right\} }\left(\mathfrak{g}\right)R_{\textrm{Par$\left(e\right)$}}\left(\mathfrak{g}\right)\in\bc\left[X_{0},\left(X_{e'}\right)_{e'\in\textrm{Anc}\left(\textrm{Par}\left(e\right)\right)}\right]\subset A\left[\Gamma\right].\]
 \end{defn}

\begin{defn}
\label{QPolyDefinition} For every $e\in\mathbf{P}\left(\Gamma\right)$,
we let \begin{eqnarray*}
Q_{e}\left(\mathfrak{g}\right) & = & S_{e}\left(\mathfrak{g}\right)R_{e}\left(\mathfrak{g}\right)\in\bc\left[X_{0},\left(X_{e'}\right)_{e'\in\textrm{Anc}\left(e\right)}\right]\subset A\left[\Gamma\right].\end{eqnarray*}

\end{defn}
\begin{defn} Given a labelled rooted tree $\mathfrak{g}=\left(\Gamma,lb\right)$,
we let \[
M\left(\mathfrak{g}\right)=\left(M_{0},\left(M_{e}\right)_{e\in\mathbf{P}\left(\Gamma\right)}\right)\in\textrm{Mat}_{2,\left(\textrm{Card}\left(\mathbf{P}\left(\Gamma\right)\right)+1\right)}\left(A\left[\Gamma\right]\right)\]
 be the matrix with the columns \[
M_{0}=\left(\begin{array}{c}
h\\
1\end{array}\right)\in\textrm{Mat}_{2,1}\left(A\left[\Gamma\right]\right)\quad\textrm{and}\quad M_{e}=\left(\begin{array}{c}
Q_{e}\left(\mathfrak{g}\right)\\
X_{e}\end{array}\right)\in\textrm{Mat}_{2,1}\left(A\left[\Gamma\right]\right),\quad e\in\mathbf{P}\left(\Gamma\right).\]
 \end{defn} 

\indent

\begin{defn}
\label{IdealDef} Given a labelled rooted tree $\mathfrak{g}=\left(\Gamma,lb\right)$,
we let $I_{\mathfrak{g}}\subset A\left[\Gamma\right]$ be the ideal
generated by the polynomials\begin{equation}
\left\{ \begin{array}{rclc}
\Delta_{0,e}\left(\mathfrak{g}\right) & = & \det\left(M_{0},M_{e}\right) & e\in\mathbf{P}\left(\Gamma\right)\\
\Delta_{e',e}\left(\mathfrak{g}\right) & = & R_{e''}^{-1}\left(\mathfrak{g}\right)\det\left(M_{e'},M_{e}\right) & \begin{array}{c}
\left(e,e'\right)\in\mathbf{P}\left(\Gamma\right)\times\textrm{Anc}\left(e\right)\\
e''=\textrm{Child}\left(e'\right)\cap\left(\downarrow e\right)\end{array}\end{array}\right..\label{FirstSimplification}\end{equation}
 The \emph{affine scheme} \emph{$V_{\mathfrak{g}}$ associated to
$\mathfrak{g}$} is the closed subscheme \[
V_{\mathfrak{g}}=\textrm{Spec}\left(B_{\mathfrak{g}}\right)\subset\textrm{Spec}\left(A\left[\Gamma\right]\right),\textrm{ with }B_{\mathfrak{g}}=A\left[\Gamma\right]/I_{\mathfrak{g}}\]

\end{defn}
\indent

\begin{rem} For every pair $\left(e,e'\right)\in\mathbf{P}\left(\Gamma\right)\times\textrm{Anc}\left(e\right)$,
\begin{eqnarray}
\Delta_{e',e}\left(\mathfrak{g}\right) & = & \left(X_{\textrm{Par}\left(e'\right)}-lb\left(e''\right)\right)X_{e}-X_{e'}Q_{e',e}\left(\mathfrak{g}\right),\label{GeneratorDef}\end{eqnarray}
 where $e''=\textrm{Child}\left(e'\right)\cap\left(\downarrow e\right)$
and $Q_{e',e}\left(\mathfrak{g}\right)=R_{f'}^{-1}\left(\mathfrak{g}\right)Q_{e}\left(\mathfrak{g}\right)\in\bc\left[\left(X_{g}\right)_{g\in\left[e',\textrm{Par}\left(e\right)\right]_{\Gamma}}\right]$. 

\end{rem}

\begin{rem} For every pair $\left(e,e'\right)\in\mathbf{P}\left(\Gamma\right)\times\textrm{Anc}\left(e\right)$
as above, the polynomials $\Delta_{0,e}\left(\mathfrak{g}\right)$,
$\Delta_{0,e'}\left(\mathfrak{g}\right)$ and $\Delta_{e',e}\left(\mathfrak{g}\right)$
satisfy the syzygy relation \begin{eqnarray}
h\Delta_{e',e}\left(\mathfrak{g}\right) & = & \left(X_{\textrm{Par}\left(e'\right)}-lb\left(e''\right)\right)\Delta_{0,e}\left(\mathfrak{g}\right)-Q_{e',e}\left(\mathfrak{g}\right)\Delta_{0,e'}\left(\mathfrak{g}\right).\label{SyzigyRelation}\end{eqnarray}
 \end{rem}

\begin{rem}
Given a pair $\left(g_{1},g_{2}\right)\in\mathbf{P}\left(\Gamma\right)\times\mathbf{P}\left(\Gamma\right)$
such that $g_{i}\not\in\left(\downarrow g_{j}\right)$, $i,j\in\left\{ 1,2\right\} $,
we let $e\in\mathbf{P}\left(\Gamma\right)$ be the first common ancestor
of $g_{1}$ and $g_{2}$ by $e$, and we let $e_{i}=\textrm{Child}\left(e\right)\cap\left(\downarrow g_{i}\right)$,
$i=1,2$. Then \begin{eqnarray}
\Delta_{g_{1},g_{2}} & := & \left(R_{e}\left(\mathfrak{g}\right)S_{e}^{\left\{ e_{1},e_{2}\right\} }\left(\mathfrak{g}\right)\right)^{-1}\det\left(M_{g_{1}},M_{g_{2}}\right)\label{SecondSimplification}\\
 & = & Q_{e,g_{1}}\left(\mathfrak{g}\right)\Delta_{e,g_{2}}\left(\mathfrak{g}\right)-Q_{e,g_{2}}\left(\mathfrak{g}\right)\Delta_{e,g_{1}}\left(\mathfrak{g}\right)\end{eqnarray}
 This means that $I_{\mathfrak{g}}$ coincides with the ideal generated
by the \emph{simplified} $2\times2$ \emph{minors} of $M\left(\mathfrak{g}\right)$,
i.e. the $2\times2$ minors of $M\left(\mathfrak{g}\right)$ simplified
by a common factor according to the rules (\ref{FirstSimplification})
and (\ref{SecondSimplification}).
\end{rem}
\begin{notation}
Hereafter, we will usually write $P$ instead of $P\left(\mathfrak{g}\right)$,
where $P\left(\mathfrak{g}\right)$ stands for any of the polynomials
$S_{e}\left(\mathfrak{g}\right)$, $R_{e}\left(\mathfrak{g}\right)$,
$Q_{e}\left(\mathfrak{g}\right)$, $\Delta_{0,e}\left(\mathfrak{g}\right)$
and $\Delta_{e',e}\left(\mathfrak{g}\right)$ attached to a labelled
rooted tree $\mathfrak{g}$. 
\end{notation}
\begin{example}
\label{Example1} Consider the following fine-labelled rooted tree

\psset{unit=1.8cm}

\begin{pspicture}(-3,-1)(3,1)

\rput(-0.5,0){$\mathfrak{g}=$}

\psline(0,0)(0.5,0.5)\psline(0,0)(1,0)\psline(0,0)(0.5,-0.5)

\psline(1,0)(1.5,0.5)\psline(1,0)(1.5,-0.5)

\rput(0,0){\textbullet}

\rput(1,0){\textbullet}\rput(0.5,0.5){\textbullet}\rput(0.5,-0.5){\textbullet}

\rput(1.5,0.5){\textbullet}\rput(1.5,-0.5){\textbullet}

\rput(-0.1,-0.2){{\small $e_0$}}

\rput(0.5,0.7){{\small $\left(e_{1,1},-1\right)$}}

\rput(0.5,-0.7){{\small $\left(e_{1,2},1\right)$}}

\rput(1.4,0){{\small $\left(e_1,0\right)$}}

\rput(1.5,0.7){{\small $\left(e_{2,1},-1\right)$}}

\rput(1.5,-0.7){{\small $\left(e_{2,2}\right)$}}

\end{pspicture}

\noindent The corresponding algebra is $A\left[\Gamma\right]=\bc\left[h\right]\left[X_{0},X_{e_{0}},X_{e_{1}}\right]$,
and the associated matrix is \[
M\left(\mathfrak{g}\right)=\left(\begin{array}{ccc}
h & X_{0}\left(X_{0}^{2}-1\right) & \left(X_{0}^{2}-1\right)\left(X_{e_{0}}^{2}-1\right)\\
1 & X_{e_{0}} & X_{e_{1}}\end{array}\right).\]
 Therefore $V_{\mathfrak{g}}$ is the Bandman and Makar-Limanov surface
\cite{BML01} with the equations \[
hX_{e_{0}}=X_{0}\left(X_{0}^{2}-1\right),\quad X_{0}X_{e_{1}}=X_{e_{0}}\left(X_{e_{0}}^{2}-1\right),\quad hX_{e_{1}}=\left(X_{0}^{2}-1\right)\left(X_{e_{0}}^{2}-1\right).\]
 This is a $GDS$ over $Z=\textrm{Spec}\left(A\right)$ for the morphism
$q_{\mathfrak{g}}:V_{\mathfrak{g}}\rightarrow Z$ induced by the inclusion
$A\hookrightarrow B_{\mathfrak{g}}$. This $\ba^{1}$-fibration $q_{\mathfrak{g}}:V_{\mathfrak{g}}\rightarrow Z$
restricts to the trivial line bundle over $Z_{*}=Z\setminus\left\{ z_{0}\right\} $,
whereas $q_{\mathfrak{g}_{1}}^{-1}\left(z_{0}\right)$ is the disjoint
union of $4$ curves \[
C_{e_{1,1}},C_{e_{1,2}}\simeq\textrm{Spec}\left(\bc\left[X_{e_{0}}\right]\right)\quad\textrm{and}\quad C_{e_{2,1}},C_{e_{2,2}}\simeq\textrm{Spec}\left(\bc\left[X_{e_{1}}\right]\right).\]
 
\end{example}

\section{$GDS$'s defined by a fine-labelled rooted trees }

\noindent This section is devoted to the proof of the following result
(see also theorem \ref{Flrt2GDSClassifTheorem} below). 

\begin{thm}
\label{EmbeddedGDSTheorem} For every \emph{fine-labelled} rooted
tree $\mathfrak{g}=\left(\Gamma,lb\right)$, the scheme $V_{\mathfrak{g}}=\textrm{Spec}\left(B_{\mathfrak{g}}\right)$
is a $GDS$ over $Z=\textrm{Spec}\left(A\right)$ for the morphism
$q_{\mathfrak{g}}:V_{\mathfrak{g}}\rightarrow Z$ induced by the inclusion
$A\hookrightarrow B_{\mathfrak{g}}$. 
\end{thm}
\begin{example}
Consider the ordinary Danielewski surface $V=V_{P,1}\subset\textrm{Spec}\left(\bc\left[x,y,z\right]\right)$
with equation $xz=P\left(y\right)$, where $P$ is a polynomial with
$r\geq1$ simple roots $y_{1},\ldots,y_{r}$. It is a $GDS$ via $pr_{x}:V\rightarrow Z=\textrm{Spec}\left(\bc\left[x\right]\right)$.
Since $xz-P\left(y\right)\in\bc\left[x,y,z\right]$ is an irreducible
polynomial, $V$ is an irreducible. By the Jacobian criterion, it
is also nonsingular. Thus $pr_{x}$ is a flat morphism. It restricts
to a trivial line bundle over $Z_{*}=\textrm{Spec}\left(\bc\left[x,x^{-1}\right]\right)$.
Since the roots of $P$ are simple, the polynomials $y-y_{i}\in\bc\left[y\right]$,
$1\leq i\leq r$ are relatively prime and so, \[
pr_{x}^{-1}\left(0\right)=\textrm{Spec}\left(\bc\left[x,y,z\right]/\left(x,xz-P\left(y\right)\right)\right)\simeq\bigsqcup_{i=1}^{r}\textrm{Spec}\left(\bc\left[y,z\right]/\left(y-y_{i}\right)\right)\]
 is the disjoint union of curves $C_{i}\simeq\textrm{Spec}\left(\bc\left[z\right]\right)$,
$1\leq i\leq r$. In general, a scheme $V_{\mathfrak{g}}$ associated
to a fine-labelled rooted tree $\mathfrak{g}$ is given by more than
one equation. Therefore, it is hard to check directly if it is irreducible
or not. In particular, a necessary condition for $q_{\mathfrak{g}}:V_{\mathfrak{g}}\rightarrow Z$
to be flat is that no irreducible or embedded component of $V_{\mathfrak{g}}$
is supported on the fiber $q_{\mathfrak{g}}^{-1}\left(z_{0}\right)$.
Let us give a direct proof of the flatness of $pr_{x}:V\rightarrow Z=\textrm{Spec}\left(\bc\left[x\right]\right)$
which illustrates the constructions used below. Since $pr_{x}$ restricts
to a trivial line bundle over $Z_{*}$, it suffices to show that $\mathcal{O}_{V,v}$
is a torsion free $\mathcal{O}_{Z,z_{0}}$-module for every point
$v\in V$ such that $pr_{x}\left(v\right)=z_{0}=Z\setminus Z_{*}$.
For every $i\in I$, we let \[
R_{i}\left(y\right)={\displaystyle \prod_{j\neq i}\left(y-y_{j}\right)=\left(y-y_{i}\right)^{-1}P\left(y\right)}\in\bc\left[x,y,z\right].\]
 Since the roots of $P\left(y\right)$ are simple, the principal open
subsets \[
\begin{array}{rcl}
U_{i}=V\cap D\left(R_{i}\right) & \simeq & \textrm{Spec}\left(\bc\left[x,y,z\right]\left[t\right]/\left(xz-P\left(y\right),R_{i}\left(y\right)t-1\right)\right)\\
 & \simeq & \textrm{Spec}\left(\bc\left[x,z,t\right]/\left(tR_{i}\left(y_{i}+xzt\right)-1\right)\right)\end{array},\quad1\leq i\leq r.\]
 cover $V$. Since $x$ does not divide $tR_{i}\left(y_{i}+xzt\right)-1$,
$U_{i}$ is flat over $Z$. Thus $V$ is covered by open subsets $U_{i}$,
$1\leq i\leq r$, flat over $Z$, whence is flat over $Z$ too. We
also deduce directly from this description that no irreducible component
of $V$ is supported on the fiber $pr_{x}^{-1}\left(0\right)$. Indeed,
the subscheme $V\left(R_{i}\right)\cap V$ is the union of the closures
of the curves \[
S_{j}=\left\{ y=y_{j},z=0,x\neq0\right\} \cap V\subset V\setminus pr_{x}^{-1}\left(0\right),\quad j\neq i,\]
 and of the irreducible components $C_{j}$, $j\neq i$, of the fiber
$pr_{x}^{-1}\left(0\right)$. This shows that for for every $1\leq i\leq r$
the generic point of $\left(C_{i}\right)_{red}$ is contained in the
closure in $V$ of the open subset $V\setminus pr_{x}^{-1}\left(0\right)$. 
\end{example}
\noindent The proof divides as follows. In subsection $1$ we show
that $q_{\mathfrak{g}}:V_{\mathfrak{g}}\rightarrow Z$ is a flat morphism.
Then is subsection $2$ we describe the fiber $q_{\mathfrak{g}}^{-1}\left(z_{0}\right)$.

\subsection*{1. Flatness of $q_{\mathfrak{g}}:V_{\mathfrak{g}}\rightarrow Z$}

\indent

\noindent In this subsection we prove the following result.

\begin{prop}
\label{TrivialFLatProp} If $\mathfrak{g}$ is a fine-labelled rooted
tree then $q_{\mathfrak{g}}:V_{\mathfrak{g}}\rightarrow Z$ is a flat
morphism restricting over $Z_{*}$ to a trivial line bundle. 
\end{prop}
\begin{enavant} \label{TrivialBundle} If $N\left(\Gamma\right)=\left\{ e_{0}\right\} $
then $B_{\mathfrak{g}}=A\left[X_{0}\right]$ and $q_{\mathfrak{g}}=pr_{Z}$
is a flat morphism. Otherwise, (\ref{SyzigyRelation}) implies that
the image of $I_{\mathfrak{g}}$ in $A_{h}\left[\Gamma\right]=A\left[\Gamma\right]\otimes_{A}A_{h}$
is generated by the polynomials\begin{equation}
\delta_{0,e}=h^{-1}\Delta_{0,e}=X_{e}-h^{-1}Q_{e},\quad e\in\mathbf{P}\left(\Gamma\right).\label{EliminationOverZ*}\end{equation}
 Since $Q_{e}$ only involves the variables $X_{0}$ and $X_{e'}$,
${\displaystyle e'\in\textrm{Anc}\left(e\right)}$, we recursively
arrive at an $A_{h}$-algebras isomorphism ${\displaystyle A_{h}\left[\Gamma\right]/I_{0}\simeq A_{h}\left[X_{0}\right]}$.
Thus $q_{\mathfrak{g}}$ restricts to a trivial line bundle over $Z_{*}\simeq\textrm{Spec}\left(A_{h}\right)$
and so, $\mathcal{O}_{V_{\mathfrak{g}},x}$ is a flat $\mathcal{O}_{Z,z}$-module
for every $x\in V_{\mathfrak{g}}$ such that $q_{\mathfrak{g}}\left(x\right)=z\in Z_{*}$.
Given $x\in V_{\mathfrak{g}}$ such that $q_{\mathfrak{g}}\left(x\right)=z_{0}$,
$\mathcal{O}_{V_{\mathfrak{g}},x}$ is a flat $\mathcal{O}_{Z,z_{0}}$-module
provided that $h$ is not a zero divisor in $\mathcal{O}_{V_{\mathfrak{g}},x}$.
So it suffices to show that $h$ is not a zero divisor in $B_{\mathfrak{g}}$.
In turn, it suffices to find $g_{1},\ldots,g_{r}\in B_{\mathfrak{g}}$
which generate the unit ideal, such that $h$ is not a zero divisor
in $\left(B_{\mathfrak{g}}\right)_{g_{i}}$ for every $1\leq i\leq r$.
The following lemma says that $V_{\mathfrak{g}}$ admits a natural
covering by principal open subsets. 

\end{enavant}

\begin{lem}
\label{CoveringLemma} The scheme $V_{\mathfrak{g}}$ is covered by
the principal open subsets \[
U_{e}=D\left(R_{e}\right)\cap V_{\mathfrak{g}}\simeq\textrm{Spec}\left(A\left[\Gamma\right]\left[T\right]/\left(I_{\mathfrak{g}},R_{e}T-1\right)\right),\quad e\in\textrm{Leaves}\left(\Gamma\right).\]

\end{lem}
\begin{proof}
Let us show that the ideal $J\subset A\left[\Gamma\right]$ generated
by the polynomials $R_{e}$, $e\in\textrm{Leaves}\left(\Gamma\right)$,
is the unit ideal of $A\left[\Gamma\right]$. Given $e'\in\mathbf{P}\left(\Gamma\right)$
such that $\textrm{Child}\left(e'\right)\subset\textrm{Leaves}\left(\Gamma\right)$,
the polynomials $S_{e'}^{\left\{ e\right\} }\in\bc\left[X_{\textrm{Par}\left(e'\right)}\right]$,
$e\in\textrm{Child}\left(e'\right)$, are relatively prime as $\mathfrak{g}$
is fine-labelled. Thus there exist polynomials $\Lambda_{e}\in A\left[\Gamma\right]$,
$e\in\textrm{Child}\left(e'\right)$, such that \[
R_{e'}=R_{e'}\sum_{e\in\textrm{Child}\left(e'\right)}\Lambda_{e}S_{e'}^{\left\{ e\right\} }={\displaystyle \sum_{e\in\textrm{Child}\left(e'\right)}\Lambda_{e}R_{e}}\]
 and so, $R_{e'}\in J$. By repeatedly applying the same argument,
we conclude that $1=R_{e_{0}}\in J$. Therefore the open subsets $D\left(R_{e}\right)$,
$e\in\textrm{Leaves}\left(\Gamma\right)$, cover $\textrm{Spec}\left(A\left[\Gamma\right]\right)$
and so, the open subsets $\left(U_{e}\right)_{e\in\textrm{Leaves}\left(\Gamma\right)}$
cover $V_{\mathfrak{g}}$. 
\end{proof}
\noindent Given $e\in\textrm{Leaves}\left(\Gamma\right)$, we have
an isomorphism of $\bc$-algebras \[
\Gamma\left(U_{e},\mathcal{O}_{V_{\mathfrak{g}}}\right)\simeq\left(B_{\mathfrak{g}}\right)_{\overline{R_{e}}}\simeq A\left[\Gamma\right]\left[T\right]/\left(I_{\mathfrak{g}},G_{e}\right),\]
 where $G_{e}=TR_{e}-1$ and $\overline{R_{e}}$ denotes the image
of $R_{e}\in A\left[\Gamma\right]$ in $B_{\mathfrak{g}}$. In view
of \ref{TrivialBundle} above, the following proposition completes
the proof of proposition \ref{TrivialFLatProp}.

\begin{prop}
\label{FlatProposition} There exists a nonzero polynomial $P_{e}\in A\left[X_{\textrm{Par}\left(e\right)},T\right]$
such that \[
\left(B_{\mathfrak{g}}\right)_{\overline{R_{e}}}\simeq A\left[X_{\textrm{Par}\left(e\right)},T\right]/\left(TP_{e}-1\right).\]
 In particular, $h$ is not a zero divisor in $\left(B_{\mathfrak{g}}\right)_{\overline{R_{e}}}$
as it does not divide $\left(TP_{e}-1\right)$. 
\end{prop}
\begin{proof}
With the notation of (\ref{MaximalChainNotation}), we let \[
L_{i}=R_{e_{i+1}}^{-1}R_{e}\in\bc\left[X_{e_{i}},\ldots,X_{e_{n_{e}-1}}\right]\subset A\left[\Gamma\right],\quad0\leq i\leq n_{e}-2.\]
 Given $e'\in\mathbf{P}\left(\Gamma\right)\setminus\left(\downarrow e\right)$,
the first common ancestor $\textrm{Anc}\left(e',e\right)$ of $e$
and $e'$ is a node $e_{i}$ for a certain indice $i\leq n_{e}-1$,
and $e''=\textrm{Child}\left(e_{i}\right)\cap\left(\downarrow e\right)\neq e_{i+1}$.
We let \[
K_{e'}=\left(X_{e_{i-1}}-lb\left(e''\right)\right)^{-1}R_{e}\in\bc\left[X_{0},X_{e_{0}},\ldots,X_{e_{n_{e}-1}}\right]\subset A\left[\Gamma\right].\]
 Modulo the ideal $I_{0}$ generated by the polynomials\[
\left\{ \begin{array}{ll}
\begin{array}{ccc}
\delta_{0,i} & = & -\left(X_{e_{i-1}}-lb\left(e_{i+1}\right)\right)+hX_{e_{i}}L_{i}T\\
 & = & TL_{i}\Delta_{0,e_{i}}+\left(X_{e_{i-1}}-lb\left(e_{i+1}\right)\right)G_{e}\end{array} & 0\leq i\leq n_{e}-2\\
\\\begin{array}{ccc}
\delta_{e'} & = & X_{e'}-TK_{e'}X_{e_{i}}Q_{e_{i},e'}\\
 & = & TK_{e'}\Delta_{e_{i},e'}-X_{e'}G_{e}\end{array} & \begin{array}{c}
e'\in\mathbf{P}\left(\Gamma\right)\setminus\left(\downarrow e\right)\\
\textrm{Anc}\left(e',e\right)=e_{i}\end{array}\end{array}\right.\]
 we can recursively eliminate all the variables in $A\left[\Gamma\right]\left[T\right]$
but $X_{\textrm{Par}\left(e\right)}$ and $T$. Therefore, there exists
a nonzero polynomial $P_{e}\in A\left[X_{\textrm{Par}\left(e\right)},T\right]$
such that $R_{e}\equiv P_{e}$ mod $I_{0}$ and so, \begin{eqnarray*}
A\left[\Gamma\right]\left[T\right]/\left(I_{0},G_{e}\right) & \simeq & A\left[X_{\textrm{Par}\left(e\right)},T\right]/\left(TP_{e}-1\right)\end{eqnarray*}
 as $G_{e}=R_{e}T-1$. In particular $h$ is not a zero divisor modulo
$\left(I_{0},G_{e}\right)$. 

\noindent Now it suffices to show that the ideals $\left(I_{0},G_{e}\right)$
and $\left(I_{\mathfrak{g}},G_{e}\right)$ of $A\left[\Gamma\right]\left[T\right]$
coincide. By construction $I_{0}\subset\left(I_{\mathfrak{g}},G_{e}\right)$.
Conversely, the identities \[
\left\{ \begin{array}{lcll}
\Delta_{0,e_{i}} & = & R_{e_{i+1}}\delta_{0,i}-hX_{e_{i}}G_{e} & 0\leq i\leq n_{e}-2\\
\Delta_{0,e'} & = & h\delta_{e'}+\left(X_{e_{i-1}}-lb\left(e_{i+1}\right)\right)^{-1}Q_{e'}\delta_{0,i} & \begin{array}{c}
e'\in\mathbf{P}\left(\Gamma\right)\setminus\left(\downarrow e\right)\\
\textrm{Anc}\left(e',e\right)=e_{i}\end{array}\end{array}\right.,\]
 guarantee that $\Delta_{0,e'}\in\left(I_{0},G_{e}\right)$ for every
$e'\in\mathbf{P}\left(\Gamma\right)$. In turn, this implies that
$h\Delta_{e'',e'}\in\left(I_{0},G_{e}\right)$ for every pair $\left(e',e''\right)\in\mathbf{P}\left(\Gamma\right)\times\textrm{Anc}\left(e'\right)$
by virtue of (\ref{SyzigyRelation}). Since $h$ is not a zero divisor
modulo $\left(I_{0},G_{e}\right)$ we conclude that $\Delta_{e'',e'}\left(\mathfrak{g}\right)\in\left(I_{0},G_{e}\right)$
for every such pair $\left(e',e''\right)$. This completes the proof. 
\end{proof}
\begin{rem}
Since $q_{\mathfrak{g}}$ is flat, (\ref{SyzigyRelation}) and (\ref{EliminationOverZ*})
implies that $B_{\mathfrak{g}}$ is a sub-$\bc$-algebra of $B_{\mathfrak{g}}\otimes_{A}A_{h}\simeq A_{h}\left[X_{0}\right]$.
On the other hand, (\ref{SyzigyRelation}) and proposition \ref{TrivialFLatProp}
assert that it suffices to add the polynomials $\Delta_{e',e}\left(\mathfrak{g}\right)$,
$\left(e,e'\right)\in\mathbf{P}\left(\Gamma\right)\times\textrm{Anc}_{\Gamma}\left(e\right)$,
to the obvious ones $\Delta_{0,e}\left(\mathfrak{g}\right)$, $e\in\mathbf{P}\left(\Gamma\right)$,
to guarantee that the fiber $q_{\mathfrak{g}}^{-1}\left(z_{0}\right)$
is the flat limit of the fibers $q_{\mathfrak{g}}^{-1}\left(z\right)$,
$z\in Z_{*}$. 
\end{rem}

\subsection*{2. The fiber $q_{\mathfrak{g}}^{-1}\left(z_{0}\right)$ }

\indent

\noindent Since $q_{\mathfrak{g}}:V_{\mathfrak{g}}\rightarrow Z$
is a flat morphism restricting to a trivial line bundle over $Z_{*}$,
$V_{\mathfrak{g}}$ is a $GDS$ over $Z$ provided that the fiber
$q_{\mathfrak{g}}^{-1}\left(z_{0}\right)$ is nonempty and reduced.
To describe this fiber, we need the following auxiliary result.

\begin{lem}
\label{FiberDecompositionLemma} If $\textrm{Child}\left(e_{0}\right)\neq\emptyset$
then there exists an isomorphism of $\bc$-algebras \[
A\left[\Gamma\right]/\left(h,I_{\mathfrak{g}}\right)\simeq\prod_{e\in\textrm{Child}\left(e_{0}\right)}\left(A\left[\Gamma\left(f\right)\right]/\left(h,I_{\mathfrak{g}\left(e\right)}\right)\right).\]

\end{lem}
\begin{enavant} The polynomial $\Delta_{0,e_{0}}\in I_{\mathfrak{g}}$
reduces to $Q_{e_{0}}$ modulo $h$. Since the roots $lb\left(e\right)\in\bc$,
$e\in\textrm{Child}\left(e_{0}\right)$, of $Q_{e_{0}}$ are simple,
we deduce that \begin{eqnarray}
A\left[\Gamma\right]/\left(h,I_{\mathfrak{g}}\right) & \simeq & \prod_{e\in\textrm{Child}\left(e_{0}\right)}\left(A\left[\Gamma\right]/\left(J_{e},I_{\mathfrak{g}}\right)\right),\label{ChineseRemainder}\end{eqnarray}
 where $J_{e}=\left(h,X_{0}-lb\left(e\right)\right)$. Thus lemma
\ref{FiberDecompositionLemma} is a consequence of the following result. 

\end{enavant} 

\begin{lem}
For every $e\in\textrm{Child}\left(e_{0}\right)$, the $\bc$-algebras
$A\left[\Gamma\right]/\left(J_{e},I_{\mathfrak{g}}\right)$ and $A\left[\Gamma\left(e\right)\right]/\left(h,I_{\mathfrak{g}\left(e\right)}\right)$
are isomorphic. 
\end{lem}
\begin{proof}
Since $X_{e_{0}}\not\in N\left(\Gamma\left(e\right)\right)$, we have
\[
A\left[\Gamma\left(e\right)\right]=A\left[X_{0}\right]\otimes_{A}\mathbf{S}\left(M_{\Gamma\left(e\right)}\right)\simeq A\left[X_{e_{0}}\right]\otimes_{A}\mathbf{S}\left(M_{\Gamma\left(e\right)}\right).\]
 With this choice of coordinates,\[
\left\{ \begin{array}{lcll}
Q_{e'}\left(\mathfrak{g}\left(e\right)\right) & = & Q_{e_{0},e'}\left(\mathfrak{g}\right) & e'\in\mathbf{P}\left(\Gamma\left(e\right)\right)\\
Q_{e'',e'}\left(\mathfrak{g}\left(e\right)\right) & = & Q_{e'',e'}\left(\mathfrak{g}\right) & \left(e',e''\right)\in\mathbf{P}\left(\Gamma\left(e\right)\right)\times\textrm{Anc}_{\Gamma\left(e\right)}\left(e'\right)\end{array}\right.\]
 We let $I_{0}\subset A\left[\Gamma\right]$ be the ideal generated
by the polynomials \[
\left\{ \begin{array}{ll}
Q_{e_{0},e'}\left(\mathfrak{g}\right) & e'\in\mathbf{P}\left(\Gamma\left(e\right)\right)\\
\Delta_{e'',e'}\left(\mathfrak{g}\right) & \left(e',e''\right)\in\mathbf{P}\left(\Gamma\left(e\right)\right)\times\left(\textrm{Anc}_{\Gamma\left(e\right)}\left(e'\right)\right)\\
\delta_{e_{0},e'}=\left(lb\left(e\right)-lb\left(g\right)\right)X_{e'}-X_{e_{0}}Q_{e_{0},e'}\left(\mathfrak{g}\right) & \left\{ \begin{array}{l}
e'\in\mathbf{P}\left(\Gamma\right)\setminus\left(\left\{ e_{0}\right\} \cup\mathbf{P}\left(\Gamma\left(e\right)\right)\right)\\
g=\textrm{Child}\left(e_{0}\right)\cap\left(\downarrow e'\right)\neq e\end{array}\right.\end{array}\right..\]
 Since $\mathfrak{g}$ is fine-labelled, $lb\left(e\right)-lb\left(g\right)\in\bc^{*}$
for every $g\in\textrm{Child}\left(e_{0}\right)\setminus\left\{ e\right\} $
and so, we can eliminate modulo $I_{0}$ all the variables $X_{e'}$
corresponding to nodes $e'\in\mathbf{P}\left(\Gamma\right)\setminus\left(\mathbf{P}\left(\Gamma\left(e\right)\right)\cup\left\{ e_{0}\right\} \right)$.
We conclude that\[
A\left[\Gamma\right]/\left(J_{e},I_{0}\right)=A\left[\Gamma\left(e\right)\right]\left[X_{0}\right]/\left(h,X_{0}-lb\left(e\right),I_{0}\right)\simeq A\left[\Gamma\left(e\right)\right]/\left(h,I_{\mathfrak{g}\left(e\right)}\right).\]
 Now it suffices to show that the ideals $\left(J_{e},I_{0}\right)$
and $\left(J_{e},I_{\mathfrak{g}}\right)$ of $A\left[\Gamma\right]$
coincide.

\noindent If $e'\in\mathbf{P}\left(\Gamma\right)\setminus\left(\mathbf{P}\left(\Gamma\left(e\right)\right)\right)$
then $\Delta_{0,e'}\left(\mathfrak{g}\right)\in J_{e}$ as $\left(X_{0}-lb\left(e\right)\right)$
divides $Q_{e'}\left(\mathfrak{g}\right)$. If $e''\in\textrm{Anc}_{\Gamma}\left(e'\right)$
then, letting $g'=\textrm{Child}\left(e''\right)\cap\left(\downarrow e'\right)$,
we obtain that \[
\Delta_{e'',e'}\left(\mathfrak{g}\right)\textrm{ mod }J_{e}=\left\{ \begin{array}{ccc}
\delta_{e_{0},e'}\in I_{0} &  & \textrm{if }e''=e_{0}\\
\left(X_{\textrm{Par}\left(e''\right)}-lb\left(g'\right)\right)\delta_{e_{0},e'}-\delta_{e_{0},e''}Q_{e'',e'}\left(\mathfrak{g}\right)\in I_{0} &  & \textrm{otherwise}\end{array}\right..\]
 If $e'\in\mathfrak{P}\left(\Gamma\left(e\right)\right)$ then $\Delta_{0,e'}\left(\mathfrak{g}\right)$
reduces to $S_{e_{0}}^{\left\{ e\right\} }\left(lb\left(e\right)\right)Q_{e_{0},e'}\left(\mathfrak{g}\right)\in I_{0}$
modulo $J_{e}$. If $e''\in\textrm{Anc}_{\Gamma}\left(e'\right)$
then either $e''=e_{0}$ and $\Delta_{e_{0},e'}\left(\mathfrak{g}\right)$
reduces to $Q_{e_{0},e'}\left(\mathfrak{g}\right)\in I_{0}$ modulo
$J_{e}$, or $e''\neq e_{0}$ and $\Delta_{e'',e'}\left(\mathfrak{g}\right)\in I_{0}$.
This proves that $\left(J_{e},I_{\mathfrak{g}}\right)=\left(J_{e},I_{0}\right)$.
\end{proof}
\noindent The following proposition completes the proof of theorem
\ref{EmbeddedGDSTheorem}. 

\begin{prop}
\label{ReducedFiberProposition} The fiber $q_{\mathfrak{g}}^{-1}\left(z_{0}\right)$
of $q_{\mathfrak{g}}:V_{\mathfrak{g}}\rightarrow Z$ is nonempty and
reduced, consisting of the disjoint union of curves $C_{e}\simeq\textrm{Spec}\left(\bc\left[X_{\textrm{Par}\left(e\right)}\right]\right)$,
$e\in\textrm{Leaves}\left(\Gamma\right)$, with defining ideals \[
I_{\mathfrak{g}}\left(e\right)=\left(I_{\mathfrak{g}},h,\left(X_{\textrm{Par}^{k+1}\left(e\right)}-lb\left(\textrm{Par}^{k-1}\left(e\right)\right)\right)_{1\leq k\leq n_{e}}\right)\subset A\left[\Gamma\right].\]
 
\end{prop}
\begin{proof}
We proceed by induction on the height $h\left(\Gamma\right)$ of $\Gamma$.
If $h\left(\Gamma\right)=0$ then $V_{\mathfrak{g}}=\textrm{Spec}\left(A\left[X_{0}\right]\right)$
and $q_{\mathfrak{g}}^{-1}\left(z_{0}\right)\simeq\textrm{Spec}\left(\bc\left[X_{0}\right]\right)$
is reduced. Otherwise, if $\textrm{Child}\left(e_{0}\right)\neq\emptyset$
then lemma \ref{FiberDecompositionLemma} implies that \[
q_{\mathfrak{g}}^{-1}\left(z_{0}\right)\simeq\textrm{Spec}\left(A\left[\Gamma\right]/\left(h,I_{\mathfrak{g}}\right)\right)\simeq\bigsqcup_{e'\in\textrm{Child}\left(e_{0}\right)}\textrm{Spec}\left(A\left[\Gamma\left(e'\right)\right]/\left(h,I_{\mathfrak{g}\left(e'\right)}\right)\right)\]
 decomposes as the disjoint union of curves $D_{e'}$ isomorphic to
the fibers $q_{\mathfrak{g}\left(e'\right)}^{-1}\left(z_{0}\right)$
of $q_{\mathfrak{g}\left(e'\right)}:V_{\mathfrak{g}\left(e'\right)}\rightarrow Z$,
$e'\in\textrm{Child}\left(e_{0}\right)$. Since $\Gamma\left(e'\right)$
has height $h\left(\Gamma\right)-1$, these fibers are nonempty and
reduced, whence $q_{\mathfrak{g}}^{-1}\left(z_{0}\right)$ is. In
view of the isomorphism (\ref{ChineseRemainder}), the description
of the irreducible components of $q_{\mathfrak{g}}^{-1}\left(z_{0}\right)$
follows easily by induction. 
\end{proof}
\begin{cor}
\label{NonsingularCoro} For every fine-labelled rooted tree $\mathfrak{g}$,
the surface $V_{\mathfrak{g}}$ is nonsingular.
\end{cor}
\begin{proof}
Indeed $q_{\mathfrak{g}}:V_{\mathfrak{g}}\rightarrow Z$ is a flat
morphism with regular geometric fibers, hence a smooth morphism. Thus
$V_{\mathfrak{g}}$ is regular as $Z\simeq\ba_{\bc}^{1}$ is. 
\end{proof}

\section{Embeddings of $GDS$'s in affine spaces}

In this section, we consider $GDS$'s $q:V\rightarrow Z=\textrm{Spec}\left(A\right)$
such that $q^{-1}\left(z\right)$ is integral for every $z\in Z_{*}$.
We prove the following result.

\begin{thm}
\label{Flrt2GDSClassifTheorem} Every $GDS$ $q:V\rightarrow Z$ is
$Z$-isomorphic to a $GDS$ $q_{\mathfrak{g}}:V_{\mathfrak{g}}\rightarrow Z$
for a certain fine-labelled rooted tree $\mathfrak{g}=\left(\Gamma,lb\right)$.
\end{thm}
\noindent The proof divides as follows. In the first subsection we
recall \cite{DubGDS} how to attach a weighted rooted tree $\gamma=\left(\Gamma,w\right)$
to a pair $\left(q:V\rightarrow Z,\phi:V\rightarrow\ba_{Z}^{1}\right)$
consisting of a $GDS$ $q:V\rightarrow Z$ and a $Z$-morphism $\phi:V\rightarrow\ba_{Z}^{1}$
restricting to an isomorphism over $Z_{*}$. In turn, this weighted
rooted tree $\gamma$ defines a $GDS$ $q_{\gamma}:W_{\gamma}\rightarrow Z$
which is $Z$-isomorphic to $V$. In the second subsection, given
a weighted rooted tree $\gamma=\left(\Gamma,w\right)$ we construct
a fine-labelled tree $\mathfrak{g}_{w}=\left(\Gamma,lb_{w}\right)$
with the same underlying tree $\Gamma$ and a closed embedding $i_{\gamma}:W_{\gamma}\hookrightarrow\textrm{Spec}\left(A\left[\Gamma\right]\right)$
inducing a $Z$-isomorphism $\psi_{\gamma}:W_{\gamma}\stackrel{\sim}{\rightarrow}V_{\mathfrak{g}_{w}}$.

\subsection{From $GDS$'s to weighted rooted trees }

\indent

\noindent In this subsection, we freely use the description of $GDS$'s
given in \cite{DubGDS}. We recall without proof how to associate
a weighted rooted tree $\gamma=\left(\Gamma,w\right)$ to a $GDS$
$q:V\rightarrow Z$. We also recall how such a weighted rooted tree
$\gamma$ defines a $GDS$ $q_{\gamma}:W_{\gamma}\rightarrow Z$ which
comes equipped with a canonical $Z$-morphism $\phi_{\gamma,0}:W_{\gamma}\rightarrow\ba_{Z}^{1}$
restricting to an isomorphism over $Z_{*}$. 

\begin{enavant} Given a $GDS$ $q:V=\textrm{Spec}\left(B\right)\rightarrow Z$,
$q$ restricts to the trivial line bundle over $Z_{*}=\textrm{Spec}\left(A_{h}\right)$.
Since $q$ is flat, we can find an isomorphism $A_{h}\left[T\right]\simeq B\otimes_{A}A_{h}$
which extends to an injection $A\left[T\right]\hookrightarrow B$.
The corresponding $Z$-morphism $\phi_{0}:V\rightarrow\ba_{Z}^{1}=\textrm{Spec}\left(A\left[T\right]\right)$
restricts to an isomorphism over $Z_{*}$. This data determines a
weighted rooted tree $\gamma=\left(\Gamma,w\right)$ as follows. 

\end{enavant}

\begin{enavant} \label{GDS2WTree} If $\phi_{0}:V\rightarrow\ba_{Z}^{1}$
is an isomorphism then we let $\gamma$ be the trivial tree with one
element $e_{0}$. Otherwise, $\phi_{0}$ is constant on the irreducible
components $C_{1},\ldots,C_{r}$ of the fiber $q^{-1}\left(z_{0}\right)$,
and the open subsets\[
V_{i}=\left(V\setminus q^{-1}\left(z_{0}\right)\right)\cup C_{i}\quad1\leq i\leq r\]
 are $Z$-isomorphic to $\ba_{Z}^{1}$. Therefore, there exist $n_{i}\geq1$,
a polynomial $\sigma_{i}\in A=\bc\left[h\right]$ of degree $\deg_{h}\left(\sigma_{i}\right)<n_{i}$
and an isomorphism $\tau_{i}:H^{0}\left(Z,q_{*}\mathcal{O}_{V_{i}}\right)\stackrel{\sim}{\rightarrow}A\left[T_{i}\right]$
such that $\phi_{0}\mid_{V_{i}}$ is induced by the polynomial $h^{n_{i}}T_{i}+\sigma_{i}\in A\left[T_{i}\right]$,
$1\leq i\leq r$. Note that over $Z_{*}$, the transition isomorphism
$\tau_{j}\circ\tau_{i}^{-1}$ is given by the $A_{h}$-algebras isomorphism\[
A_{h}\left[T_{i}\right]\stackrel{\sim}{\rightarrow}A_{h}\left[T_{j}\right],T_{i}\mapsto h^{n_{j}-n_{i}}T_{j}+h^{-n_{i}}\left(\sigma_{j}-\sigma_{i}\right)\]
 (compare with \cite{Dan89} and \cite{Fie94}). Letting \[
\sigma_{i}=\sum_{k=0}^{n_{i}-1}w_{i,k}h^{k}\in\bc\left[h\right]=A,\quad1\leq i\leq r,\]
 we consider the chains $\left(C_{i},w_{i}\right)=\left\{ e_{i,0}<e_{i,1}<\cdots<e_{i,n_{i}}\right\} $
with the weights \[
w_{i}\left(\left[e_{i,k},e_{i,k+1}\right]\right)=w_{i,k}\quad1\leq i\leq r,0\leq k\leq n_{i}-1.\]
 Since $q:V\rightarrow Z$ is affine, proposition 2.6 in \cite{DubGDS}
implies that $n_{ij}=\textrm{ord}_{z_{0}}\left(\sigma_{j}-\sigma_{i}\right)<\min\left(n_{i},n_{j}\right)$.
Consequently, there exist isomorphisms of weighted subchain\[
\theta_{ij}:C_{ij}=\left(\downarrow e_{i,n_{ij}}\right)_{C_{i}}\stackrel{\sim}{\rightarrow}C_{ji}=\left(\downarrow e_{j,n_{ij}}\right)_{C_{j}}\quad i\neq j,1\leq i,j\leq r,\]
 and a unique weighted rooted tree $\gamma=\left(\Gamma,w\right)$
with the leaves $f_{i}$, together with isomorphisms of weighted chains
$\theta_{i}:\left(C_{i},w_{i}\right)\simeq\left(\downarrow f_{i}\right)_{\Gamma}$,
$1\leq i\leq r$, such that $\theta_{i}=\theta_{j}\circ\theta_{ij}$
on $C_{ij}$ (see \cite[Proposition 1.12]{DubGDS}).

\end{enavant}

\begin{enavant} \label{WTree2GDS} Starting with a weighted rooted
tree $\gamma=\left(\Gamma,w\right)$, we construct a $GDS$ $q_{\gamma}:W_{\gamma}\rightarrow Z$
as follows. If $\Gamma$ is the trivial tree with one element $e_{0}$
then we let $W_{\gamma}=\textrm{Spec}\left(A\left[X_{0}\right]\right)\simeq\ba_{Z}^{1}$
and $q_{\gamma}=pr_{Z}$. Otherwise, we let $r\left(\gamma\right)=\textrm{Card}\left(\textrm{Leaves}\left(\Gamma\right)\right)\geq1$,
and we consider the maximal chains \begin{equation}
\left(\downarrow f\right)_{\Gamma}=\left\{ e_{f,0}=e_{0}<e_{f,1}<\cdots<e_{f,n_{f}-1}<e_{f,n_{f}}=f\right\} ,\quad f\in\textrm{Leaves}\left(\Gamma\right).\label{ChainPresentation}\end{equation}
 The prescheme $\pi_{\gamma}:X_{\gamma}\rightarrow Z$ obtained from
$Z$ by replacing the closed point $z_{0}$ by $r\left(\gamma\right)$
points $x_{f}$, $f\in\textrm{Leaves}\left(\Gamma\right)$, admits
a covering $\mathcal{U}_{\gamma}$ by the open subsets $X_{f}=\pi_{\gamma}^{-1}\left(Z_{*}\right)\cup\left\{ x_{f}\right\} \simeq Z$,
$f\in\textrm{Leaves}\left(\Gamma\right)$. Moreover $X_{ff'}=X_{f}\cap X_{f'}\simeq Z_{*}$
for every two different leaves $f$ and $f'$. We let $\mathcal{L}_{\gamma}$
be the sub-$\mathcal{O}_{X}$-module of $\mathcal{K}\left(X_{\gamma}\right)$
generated by $\left(h\circ\pi_{\gamma}\right)^{-n_{f}}$ on $X_{f}$,
$f\in\textrm{Leaves}\left(\Gamma\right)$, and we denote by $s_{\gamma}\in\Gamma\left(X,\mathcal{L}_{\gamma}\right)$
the canonical section of $\mathcal{L}_{\gamma}$ corresponding to
the constant section $1$ of $\mathcal{K}\left(X_{\gamma}\right)$.
We let $\sigma\left(\gamma\right)=\left\{ \sigma_{f}\left(\gamma\right)\right\} _{f\in\textrm{Leaves}\left(\Gamma\right)}\in C^{0}\left(\mathcal{U}_{\gamma},\mathcal{O}_{X_{\gamma}}\right)$
be the $0$-cochain defined by 

\begin{eqnarray*}
\sigma_{f}\left(\gamma\right) & = & {\displaystyle \sum_{j=0}^{n_{f}-1}w\left(\left[e_{f,j},e_{f,j+1}\right]\right)}h^{j}\in A\simeq H^{0}\left(X_{f},\mathcal{O}_{X_{\gamma}}\right),\quad f\in\textrm{Leaves}\left(\Gamma\right).\end{eqnarray*}
 Since $s_{\gamma}$ does not vanish on $X_{ff'}\simeq Z_{*}$, $f\neq f'$,
$f,f'\in\textrm{Leaves}\left(\Gamma\right)$, there exists a unique
$\check{\textrm{C}}\textrm{ech}$ $1$-cocycle $\left\{ g_{ff'}\right\} _{f,f'\in\textrm{Leaves}\left(\Gamma\right)}\in C^{1}\left(\mathcal{U}_{\gamma},\mathcal{L}_{\gamma}^{\vee}\right)$
such that \[
g_{ff'}\circ s_{\gamma}=\sigma_{f'}\left(\gamma\right)\mid_{X_{ff'}}-\sigma_{f}\left(\gamma\right)\mid_{X_{ff'}},\quad f\neq f',f,f'\in\textrm{Leaves}\left(\Gamma\right).\]

\end{enavant}

\begin{enavant} Now there exists a unique quasicoherent $\mathcal{O}_{X_{\gamma}}$-algebra
$\mathcal{A}_{\gamma}$, together with isomorphisms $\tau_{f}:\mathcal{A}_{\gamma}\mid_{X_{f}}\stackrel{\sim}{\rightarrow}\mathbf{S}\left(\mathcal{L}_{\gamma}\mid_{X_{f}}\right)$,
$f\in\textrm{Leaves}\left(\Gamma\right)$, such that over the overlaps
$X_{ff'}$, the $\mathcal{O}_{X_{ff'}}$-algebras isomorphisms $\tau_{f'}\circ\tau_{f}^{-1}:\mathbf{S}\left(\mathcal{L}_{\gamma}\mid_{X_{ff'}}\right)\stackrel{\sim}{\rightarrow}\mathbf{S}\left(\mathcal{L}_{\gamma}\mid_{X_{ff'}}\right)$
are induced by the $\mathcal{O}_{X_{ff'}}$-modules homomorphisms\begin{equation}
\left(g_{ff'},\textrm{Id}\right):\mathcal{L}_{\gamma}\mid_{X_{ff'}}\rightarrow\mathcal{O}_{X_{ff'}}\oplus\mathcal{L}_{\gamma}\mid_{X_{ff'}}\subset\mathbf{S}\left(\mathcal{L}_{\gamma}\mid_{X_{ff'}}\right).\label{BundleGlueing}\end{equation}
 By theorem 3.3 in \cite{DubGDS}, the morphism \[
q_{\gamma}=\pi_{\gamma}\circ\rho_{\gamma}:W_{\gamma}=\mathbf{Spec}\left(\mathcal{A}_{\gamma}\right)\stackrel{\rho_{\gamma}}{\rightarrow}X_{\gamma}\stackrel{\pi_{\gamma}}{\rightarrow}Z\]
 is a $GDS$ over $Z$. It restricts to the trivial line bundle over
$Z_{*}$, whereas $q_{\gamma}^{-1}\left(z_{0}\right)$ is the disjoint
union of the curves $C_{f}=\rho_{\gamma}^{-1}\left(x_{f}\right)\simeq\ba_{\bc}^{1}$,
$f\in\textrm{Leaves}\left(\Gamma\right)$. 

\end{enavant}

\begin{enavant} \label{CanonicalSection} This $GDS$ $q_{\gamma}:W_{\gamma}\rightarrow Z$
comes equipped with the canonical section $\phi_{\gamma,0}\in H^{0}\left(Z,\left(q_{\gamma}\right)_{*}\mathcal{O}_{W_{\gamma}}\right)$
whose restriction to $W_{\gamma}\mid_{X_{f}}$, $f\in\textrm{Leaves}\left(\Gamma\right)$,
is simply given by \[
\phi_{\gamma,0}\mid_{X_{f}}=\left(\sigma_{f}\left(\gamma\right),s_{\gamma}\mid_{X_{f}}\right)\in H^{0}\left(X_{f},\mathcal{A}_{\gamma}\right)\simeq H^{0}\left(X_{f},\mathbf{S}\left(\mathcal{L}_{\gamma}\right)\right).\]
 Since $s_{\gamma}$does not vanish on $\pi_{\gamma}^{-1}\left(Z_{*}\right)\subset X_{\gamma}$,
the corresponding $Z$-morphism $\phi_{\gamma,0}:W_{\gamma}\rightarrow\ba_{Z}^{1}$
restricts to an isomorphism over $Z_{*}$. If $\gamma$ is obtained
from a pair $\left(q:V\rightarrow Z,\phi:V\rightarrow\ba_{Z}^{1}\right)$
by the procedure described in \ref{GDS2WTree} then, by construction,
there exists a $Z$-isomorphism $\psi:V\stackrel{\sim}{\rightarrow}W_{\gamma}$
such that $\phi_{0}=\phi_{\gamma,0}\circ\psi$. This leads to the
following result.

\end{enavant}

\begin{thm}
\cite[Theorem 5.1]{DubGDS} For every $GDS$ $q:V\rightarrow Z$ there
exists a weighted rooted tree $\gamma=\left(\Gamma,w\right)$ and
a $Z$-isomorphism $\psi:V\stackrel{\sim}{\rightarrow}W_{\gamma}$. 
\end{thm}
\begin{example}
\label{StrangeExample1} The nonsingular surface $V\subset\textrm{Spec}\left(\bc\left[x,y,z\right]\right)$
with equation \[
x^{2}z=y^{2}-2xy-1\]
 is a $GDS$ over $Z=\textrm{Spec}\left(\bc\left[x\right]\right)$
for the morphism $pr_{x}:V\rightarrow\textrm{Spec}\left(\bc\left[x\right]\right)$.
Letting $z_{0}=\left(x\right)\in Z$, we get $V\times_{Z}Z_{*}\simeq\textrm{Spec}\left(\bc\left[x,x^{-1}\right]\left[y\right]\right)$,
whereas $pr_{x}^{-1}\left(z_{0}\right)$ is the disjoint union of
curves \[
C_{1}=\left\{ x=0,y=1\right\} \cap V\quad\textrm{and}\quad C_{2}=\left\{ x=0,y=-1\right\} \cap V.\]
 isomorphic to $\ba_{\bc}^{1}=\textrm{Spec}\left(\bc\left[z\right]\right)$.
The projection $pr_{x,y}:V\rightarrow\ba_{Z}^{1}=\textrm{Spec}\left(\bc\left[x\right]\left[y\right]\right)$
is a birational $Z$-morphism restricting to an isomorphism over $Z_{*}$.
The rational functions $T_{1}=x^{-2}\left(y-1\right)-x^{-1}$ and
$T_{2}=x^{-2}\left(y+1\right)-x^{-1}$ on $V$ induce $Z$-isomorphisms
\[
\left(V\setminus pr_{x}^{-1}\left(z_{0}\right)\right)\cup C_{1}\simeq\textrm{Spec}\left(\bc\left[x\right]\left[T_{1}\right]\right)\quad\textrm{and}\quad\left(V\setminus pr_{x}^{-1}\left(z_{0}\right)\right)\cup C_{2}\simeq\textrm{Spec}\left(\bc\left[x\right]\left[T_{2}\right]\right)\]
 Therefore, the pair $\left(pr_{x}:V\rightarrow Z,pr_{x,y}:V\rightarrow\ba_{Z}^{1}\right)$
corresponds to the following weighted rooted tree 

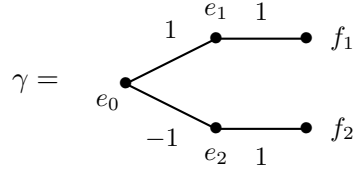
\psset{unit=1.2cm}

\begin{pspicture}(-4.5,-1.5)(6,1.2)

\rput(-1,0){$\gamma=$}

\psline(0,0)(1,0.5)\psline(0,0)(1,-0.5)

\psline(1,0.5)(2,0.5)\psline(1,-0.5)(2,-0.5)

\rput(0,0){\textbullet}

\rput(1,0.5){\textbullet}\rput(2,0.5){\textbullet}

\rput(1,-0.5){\textbullet}\rput(2,-0.5){\textbullet}

\rput(-0.2,-0.2){{\small $e_0$}}

\rput(1,0.8){{\small $e_1$}}

\rput(1,-0.8){{\small $e_2$}}

\rput(2.4,-0.5){{\small $f_2$}}

\rput(2.4,0.5){{\small $f_1$}}

\rput(0.5,0.6){{\small $1$}}

\rput(0.4,-0.6){{\small $-1$}}

\rput(1.5,0.8){{\small $1$}}

\rput(1.5,-0.8){{\small $1$}}

\end{pspicture}
\end{example}

\subsection{Embedding of a $GDS$ $W_{\gamma}$ in an affine space }

\indent

\noindent Given a $GDS$ $q_{\gamma}:W_{\gamma}\rightarrow Z$ defined
by a weighted rooted tree $\gamma=\left(\Gamma,w\right)$, we construct
in \ref{PropertyofPhi_0}-\ref{EMbeddingTheorem} below a fine-labelled
tree $\mathfrak{g}=\left(\Gamma,lb_{w}\right)$ with the same underlying
rooted tree $\Gamma$, and a collection of regular functions \[
\left(\phi_{\gamma,0},\left(\phi_{\gamma,e}\right)_{e\in\mathbf{P}\left(\Gamma\right)}\right)\in B_{\gamma}=H^{0}\left(Z,\left(q_{\gamma}\right)_{*}\mathcal{O}_{W_{\gamma}}\right)\]
 defining a closed embedding $i_{\gamma}:W_{\gamma}\hookrightarrow\ba_{Z}^{d\left(\Gamma\right)+1}=\textrm{Spec}\left(A\left[\Gamma\right]\right)$
with image $V_{\mathfrak{g}_{w}}$. \begin{enavant}\label{PropertyofPhi_0}
Consider the canonical section $\phi_{\gamma,0}\in B_{\gamma}$ defined
in \ref{CanonicalSection}. If $\Gamma$ is the trivial tree with
one element $e_{0}$ then \[
W_{\gamma}\simeq\textrm{Spec}\left(A\left[\phi_{\gamma,0}\right]\right)\simeq\textrm{Spec}\left(A\left[\Gamma\right]\right).\]
 Otherwise, for every $f\in\textrm{Leaves}\left(\Gamma\right)$, the
multiplication by $h^{n_{f}}$ gives rise to an $\mathcal{O}_{X_{f}}$-modules
isomorphism $\mathcal{L}_{\gamma}\mid_{X_{f}}\stackrel{\sim}{\rightarrow}\mathcal{O}_{X_{f}}$,
whence to an $A$-algebras isomorphism $\tau_{f}:H^{0}\left(X_{f},\mathcal{A}_{\gamma}\right)\stackrel{\sim}{\longrightarrow}A\left[T\right]$
such that \begin{eqnarray}
\tau_{f}\left(\phi_{\gamma,0}\right) & = & h^{n_{f}}T+\sigma_{f}\left(\gamma\right)=\sum_{k=0}^{n_{f}}w_{f,k}h^{k}\in A\left[T\right],\label{Intrivilization}\end{eqnarray}
 where $w_{f,n_{f}}=T$. Since $n_{f}\geq1$, we deduce that $\phi_{\gamma,0}$
is locally constant on $q_{\gamma}^{-1}\left(z_{0}\right)$, with
the value $w_{f,0}\in\bc$ on $C_{f}=\rho_{\gamma}^{-1}\left(x_{f}\right)\simeq\textrm{Spec}\left(\bc\left[T\right]\right)$.
For every $e\in N\left(\Gamma\right)$, we let \[
F\left(e\right)={\displaystyle \bigcup_{f\in\textrm{Leaves}\left(\Gamma\left(e\right)\right)}C_{f}\subset q_{\gamma}^{-1}\left(z_{0}\right)}.\]
 In particular, $F\left(e_{0}\right)=q_{\gamma}^{-1}\left(z_{0}\right)$
and $F\left(f\right)=C_{f}$ for every $f\in\textrm{Leaves}\left(\Gamma\right)$.
We have the following result. 

\end{enavant}

\begin{prop}
\label{EmbeddingMainProp} If $h\left(\Gamma\right)\geq1$ then there
exists a fine-labelled rooted tree $\mathfrak{g}_{w}=\left(\Gamma,lb_{w}\right)$
with underlying tree $\Gamma$, and a collection of regular functions
$\phi_{\gamma,\textrm{Par}\left(e_{0}\right)}=\phi_{\gamma,0}$,$\phi_{\gamma,e}\in B_{\gamma}$,
$e\in\mathbf{P}\left(\Gamma\right)$, such that the following hold.

1) For every $e\in\mathbf{P}\left(\Gamma\right)$ with $\left(\downarrow e\right)_{\Gamma}=\left\{ e_{0}<e_{1}<\cdots<e_{n}=e\right\} $, 

\[
\phi_{\gamma,e}=h^{-1}Q_{e}\left(\mathfrak{g}_{w}\right)\left(\phi_{\gamma,0},\phi_{\gamma,e_{0}},\ldots,\phi_{\gamma,e_{n-1}}\right).\]

2) Letting $\textrm{Child}\left(e\right)=\left\{ g_{1},\ldots,g_{r}\right\} $,
$\phi_{\gamma,\textrm{Par}\left(e\right)}$ is constant on $F\left(g_{i}\right)$
with the value $lb_{w}\left(g_{i}\right)\in\bc$ whereas $\phi_{\gamma,e}$
restricts to a coordinate function on $C_{g}\subset q_{\gamma}^{-1}\left(z_{0}\right)$
for every $g\in\textrm{Child}\left(e\right)\cap\textrm{Leaves}\left(\Gamma\right)$. 
\end{prop}
\begin{enavant} Starting with $\phi_{\gamma,0}$, we construct by
induction the labelling $lb_{w}:N\left(\Gamma\right)\setminus\left\{ e_{0}\right\} \rightarrow\bc$
and the regular functions $\phi_{\gamma,e}$, $e\in\mathbf{P}\left(\Gamma\right)$.
More precisely, for every $n\leq h\left(\Gamma\right)$ in $N\left(\Gamma\right)$,
we define the labelling function $lb_{w}$ on $N_{n}\left(\Gamma\right)$
in terms of the functions $\phi_{\gamma,e'}$ corresponding to nodes
$e'\in N_{k}\left(\Gamma\right)$, $k\leq n-2$. Then we define the
functions $\phi_{\gamma,e}$ for nodes $e\in N_{n-1}\left(\Gamma\right)$
in terms of those functions $\phi_{\gamma,e'}$ and the values of
$lb_{w}$ on the nodes $e''\in N_{k}\left(\Gamma\right)$, $k\leq n$. 

\end{enavant}

\begin{enavant}\label{Step1} Step 1. Given $g\in\textrm{Child}\left(e_{0}\right)$,
we let \[
lb_{w}\left(g\right):=\phi_{\gamma,0}\mid_{F\left(g\right)}=w_{f,0},\]
 where $f$ is any leaf of $\Gamma\left(g\right)$. This is well-defined
as $g=e_{f,1}$ for every such leaf $f$ (see (\ref{ChainPresentation})).
If $g'$ is another child of $e_{0}$ then there exists $f'\in\textrm{Leaves}\left(\Gamma\left(g'\right)\right)$
such that $g'=e_{f',1}$, and hence \[
lb_{w}\left(g\right)=w_{f,0}\neq w_{f',0}=lb_{w}\left(g'\right)\]
 by definition \ref{WeightDef}. This defines the labelling function
$lb_{w}$ on the nodes at level $1$. Clearly, \begin{eqnarray*}
\phi_{\gamma,e_{0}}: & = & h^{-1}Q_{e_{0}}\left(\mathfrak{g}_{w}\right)\left(\phi_{\gamma,0}\right)=h^{-1}\prod_{g\in\textrm{Child}\left(e_{0}\right)}\left(\phi_{\gamma,0}-lb_{w}\left(g\right)\right)\in B_{\gamma}\otimes_{A}A_{h}\end{eqnarray*}
 is a regular function on $W_{\gamma}$. Letting \[
s_{g}=S_{e_{0}}^{\left\{ g\right\} }\left(\mathfrak{g}_{w}\right)\left(lb_{w}\left(g\right)\right)=\prod_{g'\in\textrm{Child}\left(e_{0}\right)\setminus\left\{ g\right\} }\left(lb_{w}\left(g\right)-lb_{w}\left(g'\right)\right)\in\bc^{*},\quad g\in\textrm{Child}\left(e_{0}\right),\]
we deduce from Taylor's Formula that for every $f\in\textrm{Leaves}\left(\Gamma\left(g\right)\right)$
there exists $P_{0}\left(g,f\right)\in A\left[T\right]$ such that
\[
\tau_{f}\left(\phi_{\gamma,e_{0}}\right)=s_{g}w_{f,1}+hP_{0}\left(g,f\right)\in A\left[T\right].\]
 If $g$ is a leaf of $\Gamma$ then $w_{g,1}=T$ (see (\ref{Intrivilization}))
and so, $\phi_{\gamma,e_{0}}$ restricts to a coordinate function
on $C_{g}\subset q_{\gamma}^{-1}\left(z_{0}\right)$. Thus (1)-(2)
hold for $e_{0}\in N\left(\Gamma\right)$.

\end{enavant}

\begin{enavant} \label{Step2} Step 2. To define $lb_{w}$ on the
nodes at level $2$, we proceed as follows. Given $g'\in\textrm{Child}\left(g\right)$,
$\phi_{\gamma,e_{0}}$ is constant on $F\left(g'\right)\subset q_{\gamma}^{-1}\left(z_{0}\right)$
with the value \[
lb_{w}\left(g'\right):=s_{g}w_{f,1}\in\bc,\]
 where $f$ is any leaf of $\Gamma\left(g'\right)$. Given another
child $g''$ of $g$ and a leaf $f'$ of $\Gamma\left(g''\right)$,
$g=e_{f,1}=e_{f',1}$ is the first common ancestor of $f$ and $f'$.
Thus $e_{f,2}\neq e_{f',2}$, whence $w_{f,1}\neq w_{f',1}$ as $\gamma$
is a weighted rooted tree. Since $s_{g}\in\bc^{*}$ we deduce that
$lb_{w}\left(g'\right)\neq lb_{w}\left(g''\right)$. This defines
the labelling function $lb_{w}$ on the nodes at level $2$. 

\end{enavant}

\noindent Suppose that the labelling $lb_{w}$ has been defined for
every node $e$ at level $k\leq n+1$ and that for every $e\in\mathbf{P}\left(\Gamma\right)\cap N_{n}\left(\Gamma\right)$
with $\left(\downarrow e\right)=\left\{ e_{0}<e_{1}<\cdots<e_{n}=e\right\} $,
the regular functions $\phi_{e_{i}}\in B_{\gamma}$, $0\leq i\leq n-1$
have been constructed, satisfying (1) and (2) of proposition \ref{EmbeddingMainProp}.
Since $lb_{w}$ is defined on $N_{k}\left(\Gamma\right)$, $k\leq n+1$,
the polynomial $Q_{e}\left(\mathfrak{g}_{w}\right)\in\bc\left[X_{0},X_{e_{0}},\ldots,X_{e_{n-1}}\right]$
is well-defined (see \ref{QPolyDefinition}). We have the following
result. 

\begin{lem}
If $e$ is not a leaf of $\Gamma$ then \[
\phi_{\gamma,e_{n}}=h^{-1}Q_{e}\left(\mathfrak{g}_{w}\right)\left(\phi_{\gamma,0},\phi_{\gamma,e_{0}},\ldots,\phi_{\gamma,e_{n-1}}\right)\in B_{\gamma}\otimes_{A}A_{h}\]
 is a regular function on $W_{\gamma}$ satisfying (1)-(2) in proposition
\ref{EmbeddingMainProp}. 
\end{lem}
\begin{proof}
In view of the definition of the functions $\phi_{\gamma,e_{i}}$,
$0\leq i\leq n-1$, we deduce from Taylor's Formula that for every
$g\in\textrm{Child}\left(e\right)$ and every $f\in\textrm{Leaves}\left(\Gamma\left(g\right)\right)$,
there exists $P_{n-1}\left(g,f\right)\in A\left[T\right]$ such that
\[
\tau_{f}\left(\phi_{\gamma,e_{n-1}}\right)=lb_{w}\left(g\right)+\left(\prod_{i=0}^{n-1}s_{e_{i}}^{n-1-i}w_{f,n+1}+Q_{n-1,g}\left(w_{f,0},\ldots,w_{f,n}\right)\right)h+h^{2}P_{n-1}\left(g,f\right)\in A\left[T\right],\]
 where $Q_{n-1,g}\in\bc\left[Y_{0},\ldots,Y_{n}\right]$ is independent
on the choice of a leaf $f$ of $\Gamma\left(g\right)$ and \[
s_{e_{i}}=S_{e_{i}}^{\left\{ e_{i+1}\right\} }\left(\mathfrak{g}_{w}\right)\left(lb_{w}\left(e_{i+1}\right)\right)\in\bc^{*},\quad1\leq i\leq n-1.\]
 By definition,\[
\phi_{\gamma,e_{n}}=h^{-1}S_{e_{n}}\left(\mathfrak{g}_{w}\right)\left(\phi_{\gamma,e_{n-1}}\right)\prod_{i=0}^{n-1}S_{e_{i}}^{\left\{ e_{i+1}\right\} }\left(\mathfrak{g}_{w}\right)\left(\phi_{\gamma,e_{i-1}}\right)\in B_{\gamma}\otimes_{A}A_{h}\]
 and $\phi_{\gamma,e_{i-1}}\mid_{C\left(f\right)}=lb_{w}\left(e_{i+1}\right)$,
$0\le i\leq n$. Again, Taylor's Formula implies that there exists
$P_{n}\left(g,f\right)\in A\left[T\right]$ such that \begin{equation}
\tau_{f}\left(\phi_{\gamma,e_{n}}\right)=\left(\lambda_{g}w_{f,n+1}+\mu_{g}\right)+hP_{n}\left(g,f\right)\in A\left[T\right],\label{TaylorDev}\end{equation}
 where \[
\left\{ \begin{array}{lcl}
\lambda_{g} & = & S_{e_{n}}^{\left\{ g\right\} }\left(\mathfrak{g}_{w}\right)\left(lb_{w}\left(g\right)\right)\left({\displaystyle \prod_{i=0}^{n-1}}s_{e_{i}}^{n-i}\right)\in\bc^{*},\\
\mu_{g} & = & S_{e_{n}}^{\left\{ g\right\} }\left(\mathfrak{g}_{w}\right)\left(lb_{w}\left(g\right)\right)\left({\displaystyle \prod_{i=0}^{n-1}}s_{e_{i}}\right)Q_{n-1,g}\left(w_{f,0},\ldots,w_{f,n}\right)\in\bc.\end{array}\right.\]
Thus $\phi_{\gamma,e_{n}}$ is regular on $W_{\gamma}\mid_{X_{f}}$
for every $f\in\textrm{Leaves}\left(\Gamma\left(e_{n}\right)\right)$.
Given $f'\in\textrm{Leaves}\left(\Gamma\right)\setminus\textrm{Leaves}\left(\Gamma\left(e_{n}\right)\right)$,
there exists an indice $k\leq n-1$ such that $e_{k}$ is the first
common ancestor of $f'$ and $e_{n}$. Letting $e'=\textrm{Child}\left(e_{k}\right)\cap\left(\downarrow f'\right)$,
$\left(X_{e_{k-1}}-lb_{w}\left(e'\right)\right)$ divides $Q_{e_{n}}\left(\mathfrak{g}_{w}\right)$
and so, $\tau_{f'}\left(\phi_{\gamma,e_{k-1}}-lb_{w}\left(e'\right)\right)\in hA\left[T\right]$.
In turn, this implies that $\tau_{f'}\left(\phi_{\gamma,e_{n}}\right)\in A\left[T\right]$.
Thus $\phi_{\gamma,e_{n}}$ is a regular function on $W_{\gamma}$.
If $g\in\textrm{Child}\left(e_{n}\right)$ is a leaf of $\Gamma$
then $w_{g,n+1}=T$ (see (\ref{Intrivilization})) and so, $\phi_{\gamma,e_{n}}$
restricts to a coordinate function on $C_{g}\subset q_{\gamma}^{-1}\left(z_{0}\right)$. 
\end{proof}
\begin{enavant} Consider the regular function $\phi_{\gamma,e_{n}}$
constructed above. Given $g\in\textrm{Child}\left(e_{n}\right)\setminus\textrm{Leaves}\left(\Gamma\right)$
and $g'\in\textrm{Child}\left(g\right)$, we deduce from (\ref{TaylorDev})
that $\phi_{\gamma,e_{n}}$ is constant on $F\left(g'\right)\subset q_{\gamma}^{-1}\left(z_{0}\right)$
with the value \begin{equation}
lb_{w}\left(g'\right):=\lambda_{g}w_{f,n+1}+\mu_{g},\label{LabelAt(n+2)}\end{equation}
 where $f$ is any leaf of $\Gamma\left(g'\right)$. Since neither
$\lambda_{g}$ nor $\mu_{g}$ depend on the choice of $f\in\textrm{Leaves}\left(\Gamma\left(g\right)\right)$,
the same argument as in \ref{Step1} shows that $lb_{w}\left(g'\right)\neq lb_{w}\left(g''\right)$
for every two distinct children $g'$ and $g''$ of $g$. In this
way, we define the labelling function $lb_{w}$ on the nodes of $\Gamma$
at level $n+2$. Now the proof of proposition \ref{EmbeddingMainProp}
can be completed by induction. 

\end{enavant}

\begin{rem}
\label{InversionRemark} By construction, the labelling $lb_{w}:N\left(\Gamma\right)\rightarrow\bc$
is uniquely determined by the weight function $w:E\left(\Gamma\right)\rightarrow\bc$.
Conversely, since $lb_{w}\left(e\right)=w\left(\left[e_{0},e\right]_{\Gamma}\right)$
for every node $e\in N_{1}\left(\Gamma\right)$ (see \ref{Step1}),
we deduce from (\ref{LabelAt(n+2)}) that the weight function $w:E\left(\Gamma\right)\rightarrow\bc$
can be recursively recovered from the associated labelling $lb_{w}$. 
\end{rem}
\noindent The following result completes the proof of theorem \ref{Flrt2GDSClassifTheorem}.

\begin{thm}
\label{EMbeddingTheorem} For every $GDS$ $q_{\gamma}:W_{\gamma}\rightarrow Z=\textrm{Spec}\left(A\right)$
defined from a weighted rooted tree $\gamma=\left(\Gamma,w\right)$,
the $A$-algebras homomorphism \[
A\left[\Gamma\right]\rightarrow B_{\gamma},\left\{ \begin{array}{ll}
X_{0}\mapsto\phi_{\gamma,0}\\
X_{e}\mapsto\phi_{\gamma,e} & e\in\mathbf{P}\left(\Gamma\right)\end{array}\right.\]
 is surjective, with kernel $I_{\mathfrak{g}_{w}}$. In other words,
the functions $\phi_{\gamma,0},\left(\phi_{\gamma,e}\right)_{e\in\mathbf{P}\left(\Gamma\right)}$
define a closed embedding $i_{\gamma}:W_{\gamma}\hookrightarrow\ba_{Z}^{d\left(\Gamma\right)+1}=\textrm{Spec}\left(A\left[\Gamma\right]\right)$
inducing a $Z$-isomorphism $\psi_{\gamma}:W_{\gamma}\stackrel{\sim}{\rightarrow}V_{\mathfrak{g}_{w}}$. 
\end{thm}
\begin{proof}
By construction, $W_{\gamma}\times_{Z}Z_{*}\simeq\textrm{Spec}\left(A_{h}\left[\phi_{\gamma,0}\right]\right)$.
Given two different leaves $f$ and $f'$ of $\Gamma$ with first
common ancestor $e\in\mathbf{P}\left(\Gamma\right)$, it follows from
proposition \ref{EmbeddingMainProp} that $\phi_{\gamma,\textrm{Par}\left(e\right)}$
takes disctint values on the curves $C_{f}$ and $C_{f'}$. Thus $i_{\gamma}$
separates the irreducible components of $q_{\gamma}^{-1}\left(z_{0}\right)$.
If $f\in\textrm{Leaves}\left(\Gamma\right)$ then $\phi_{\gamma,\textrm{Par}\left(f\right)}$
restricts to a coordinate function on $C_{f}$. This proves that $i_{\gamma}$
is a closed embedding of $Z$-schemes. By (1) in proposition \ref{EmbeddingMainProp},
the relation \[
\Delta_{0,e}\left(\mathfrak{g}_{w}\right)\left(\phi_{\gamma,0},\left(\phi_{\gamma,e}\right)_{e\in\mathfrak{P}\left(\Gamma\right)}\right)=h\phi_{\gamma,e}-Q_{e}\left(\mathfrak{g}_{w}\right)\left(\phi_{\gamma,0},\left(\phi_{\gamma,e}\right)_{e\in\mathfrak{P}\left(\Gamma\right)}\right)=0\]
 holds in $B_{\gamma}$ for every $e\in\mathbf{P}\left(\Gamma\right)$.
In turn, (\ref{SyzigyRelation}) implies that \[
h\Delta_{e',e}\left(\mathfrak{g}_{w}\right)\left(\phi_{\gamma,0},\left(\phi_{\gamma,e}\right)_{e\in\mathfrak{P}\left(\Gamma\right)}\right)=0\quad\left(e,e'\right)\in\mathbf{P}\left(\Gamma\right)\times\textrm{Anc}_{\Gamma}\left(e\right)\]
 and so, $\Delta_{e',e}\left(\mathfrak{g}_{w}\right)\left(\phi_{0},\left(\phi_{e}\right)_{e\in\mathfrak{P}\left(\Gamma\right)}\right)=0$
as $B_{\gamma}$ is an integral $A$-algebra. Since $i_{\gamma}$
maps the irreducible components of $q_{\gamma}^{-1}\left(z_{0}\right)$
bijectively on the irreducible components of $q_{\mathfrak{g}_{w}}^{-1}\left(z_{0}\right)$
(see proposition \ref{ReducedFiberProposition}), we deduce that $i_{\gamma}$
induces a bijective morphism $\psi_{\gamma}:W_{\gamma}\rightarrow V_{\mathfrak{g}_{w}}$.
By virtue of Zariski's Main Theorem, $\psi_{\gamma}$ is an isomorphism
as $V_{\mathfrak{g}_{w}}$ is nonsingular (see corollary \ref{NonsingularCoro}),
whence normal. 
\end{proof}
\begin{example}
We again consider the surface $V\subset\textrm{Spec}\left(\bc\left[x,y,z\right]\right)$
with equation \[
x^{2}z=y^{2}-2xy-1\]
 introduced in example \ref{StrangeExample1}. Starting with the weighted
rooted tree $\gamma=\left(\Gamma,w\right)$ corresponding to the pair
\[
\left(pr_{x}:V\rightarrow Z,pr_{x,y}:V\rightarrow\ba_{Z}^{1}\right),\]
 the above procedure gives the following fine-labelled rooted tree 

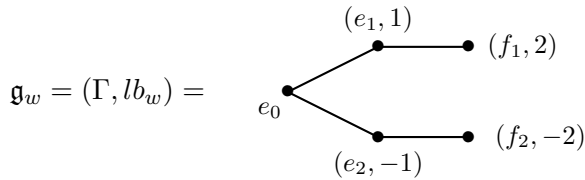
\psset{unit=1.2cm}

\begin{pspicture}(-4.5,-1)(6,1.2)

\rput(-2,0){$\mathfrak{g}_w=\left(\Gamma,lb_w\right)= $}

\psline(0,0)(1,0.5)\psline(0,0)(1,-0.5)

\psline(1,0.5)(2,0.5)\psline(1,-0.5)(2,-0.5)

\rput(0,0){\textbullet}

\rput(1,0.5){\textbullet}\rput(2,0.5){\textbullet}

\rput(1,-0.5){\textbullet}\rput(2,-0.5){\textbullet}

\rput(-0.2,-0.2){{\small $e_0$}}

\rput(1,0.8){{\small $\left(e_1,1\right)$}}

\rput(1,-0.8){{\small $\left(e_2,-1\right)$}}

\rput(2.8,-0.5){{\small $\left(f_2,-2\right)$}}

\rput(2.6,0.5){{\small $\left(f_1,2\right)$}}

\end{pspicture}
\end{example}
\noindent Thus $V_{\mathfrak{g}_{w}}\subset\textrm{Spec}\left(A\left[X_{0}\right]\left[X_{e_{0}},X_{e_{1}},X_{e_{2}}\right]\right)$
is the surface with equations\[
\left\{ \begin{array}{c}
hX_{e_{0}}=\left(X_{0}^{2}-1\right),\quad hX_{e_{1}}=\left(X_{0}+1\right)\left(X_{e_{0}}-2\right),\quad hX_{e_{2}}=\left(X_{0}-1\right)\left(X_{e_{0}}+2\right)\\
\left(X_{0}-1\right)X_{e_{1}}=X_{e_{0}}\left(X_{e_{0}}-2\right),\quad\left(X_{0}+1\right)X_{e_{2}}=X_{e_{0}}\left(X_{e_{0}}+2\right)\end{array}\right..\]
 This shows that in general, the embedding of a $GDS$ $V$ obtained
by the above procedure is not ''the best possible''. However, note
that the $\bc$-algebras homomorphism \[
A\left[X_{0},X_{e_{0}},X_{e_{1}},X_{e_{2}}\right]\rightarrow\bc\left[x,x^{-1}y,z\right],\left\{ \begin{array}{l}
h\mapsto x,\quad X_{0}\mapsto y\\
X_{e_{0}}\mapsto xz+2y\\
X_{e_{1}}\mapsto x^{-1}\left(xz+2y-2\right)\left(y+1\right)\\
X_{e_{2}}\mapsto x^{-1}\left(xz+2y+2\right)\left(y-1\right)\end{array}\right.\]
 induces an isomorphism $V\stackrel{\sim}{\rightarrow}V_{\mathfrak{g}}$. 

\indent

\noindent The following result is complementary to proposition \ref{EmbeddingMainProp}.

\begin{prop}
\label{Flrt2FwrtTheorem} For every $GDS$ $q_{\mathfrak{g}}:V_{\mathfrak{g}}\rightarrow Z$
defined by a fine-labelled rooted tree $\mathfrak{g}=\left(\Gamma,lb\right)$,
the weighted tree $\gamma$ associated to the pair $\left(q_{\mathfrak{g}}:V_{\mathfrak{g}}\rightarrow Z,pr_{X_{0}}:V_{\mathfrak{g}}\rightarrow\ba_{Z}^{1}\right)$
has the same underlying tree $\Gamma$ as $\mathfrak{g}$, and the
$GDS$ $V_{\mathfrak{g}}$ is $Z$-isomorphic to $W_{\gamma}$ via
the closed embedding $i_{\gamma}:W_{\gamma}\hookrightarrow\ba_{Z}^{d\left(\Gamma\right)+1}$.
\end{prop}
\begin{proof}
By construction of the embedding $i_{\gamma}:W_{\gamma}\hookrightarrow\ba_{Z}^{d\left(\Gamma\right)+1}$,
the canonical $Z$-morphism $\phi_{0,\gamma}:W_{\gamma}\rightarrow\ba_{Z}^{1}$
factors as $\phi_{\gamma,0}=pr_{X_{0}}\circ i_{\gamma}$ (see \ref{EMbeddingTheorem}).
So the statement follows from remark \ref{InversionRemark}. 
\end{proof}

\section{$\gaz$-actions on $GDS$'s over $Z$}

In this subsection we study $\gaz$-actions on a $GDS$ $q:V\rightarrow Z$.
In view of the well known correspondence between algebraic $\gaz$-actions
on $V$ and locally nilpotent $A$-derivations of the algebra $B=H^{0}\left(Z,q_{*}\mathcal{O}_{V}\right)$,
this is the same as to study the set $\textrm{LND}_{A}\left(B\right)$
of these locally nilpotent derivations.

\subsection{$\gaz$-actions on the $GDS$'s $q_{\gamma}:W_{\gamma}\rightarrow Z$ }

\indent

\noindent Here we recall following \cite{DubGDS} the description
of $\gaz$-actions on a $GDS$ $q_{\gamma}:W_{\gamma}\rightarrow Z$
defined from a weighted rooted tree $\gamma=\left(\Gamma,w\right)$.
By construction (see \ref{WTree2GDS}), the morphism $q_{\gamma}:W_{\gamma}\rightarrow Z$
factors through $\rho_{\gamma}:W_{\gamma}\rightarrow X_{\gamma}$,
in such a way that $W_{\gamma}\mid_{X_{f}}$ is isomorphic to the
restriction of the line bundle $p:L_{\gamma}=\mathbf{Spec}\left(\mathbf{S}\left(\mathcal{L}_{\gamma}\right)\right)\rightarrow X_{\gamma}$
over $X_{f}$, $f\in\textrm{Leaves}$. This line bundle $L_{\gamma}$
has a natural structure of group scheme as it represents the group
functor \[
\left(Sch_{/X_{\gamma}}\right)\rightarrow\left(Grp\right),\left(Y\stackrel{f}{\rightarrow}X_{\gamma}\right)\mapsto H^{0}\left(Y,f^{*}\mathcal{L}_{\gamma}^{\vee}\right),\]
 and $W_{\gamma}$ comes naturally equipped with a left action $\xi_{\gamma}:L_{\gamma}\times_{X_{\gamma}}W_{\gamma}\rightarrow W_{\gamma}$
making $\rho_{\gamma}:W_{\gamma}\rightarrow X_{\gamma}$ a principal
homogeneous $L_{\gamma}$-bundle. 

\begin{enavant} \label{Function2Action}Every nonzero section $s\in H^{0}\left(X_{\gamma},\mathcal{L}_{\gamma}^{\vee}\right)$
gives rise to a group homomorphism \[
\phi_{s}:\mathbb{G}_{a,X_{\gamma}}=\mathbf{Spec}\left(\mathcal{O}_{X_{\gamma}}\left[T\right]\right)\rightarrow L_{\gamma},\]
 whence to a nontrivial $\mathbb{G}_{a,X_{\gamma}}$-action $\mu_{s}=\xi_{\gamma}\circ\left(\phi_{s}\times\textrm{Id}\right)$
on $W_{\gamma}$. By proposition 3.10 in \cite{DubGDS}, every nontrivial
$\gaz$-action on $q_{\gamma}:W_{\gamma}\rightarrow Z$ is induced
by such a nonzero section $s$. Since $K\left(X_{\gamma}\right)\simeq K\left(Z\right)\simeq\textrm{Frac}\left(A\right)$,
\[
H^{0}\left(X_{\gamma},\mathcal{L}_{\gamma}^{\vee}\right)\simeq\left\{ g\in K\left(X_{\gamma}\right),\textrm{div}\left(g\right)-\sum_{f\in\textrm{Leaves}\left(\Gamma\right)}n_{f}x_{f}\geq0\right\} \cup\left\{ 0\right\} \]
 and so, nonzero sections $s\in H^{0}\left(X_{\gamma},\mathcal{L}_{\gamma}^{\vee}\right)$
are in one-to-one correspondence with regular functions $h^{m}g\in A$,
where $m\geq h\left(\Gamma\right)=\max\left(n_{f}\right)_{f\in\textrm{Leaves}\left(\Gamma\right)}$
and $g\in A\setminus hA$. 

\end{enavant}

\begin{rem}
\label{Co-actionRem} Let $s\in H^{0}\left(X_{\gamma},\mathcal{L}_{\gamma}^{\vee}\right)$
be a nonzero section corresponding to a regular function $h^{m}g\in A$.
Over $X_{f}$, the multiplication by $h^{n_{f}}$ induces an isomorphism
of $\mathcal{O}_{X_{f}}$-algebras $\tau_{f}:\mathcal{A}_{\gamma}\mid_{X_{f}}\stackrel{\sim}{\rightarrow}\mathcal{O}_{X_{f}}\left[U\right]$
such that the group co-action corresponding to the restriction of
the $\gaz$-action $\mu_{s}$ on $W_{\gamma}\mid_{X_{f}}$ is given
by the $\mathcal{O}_{X_{f}}$-algebras homomorphism \[
\mathcal{O}_{X_{f}}\left[U\right]\rightarrow\mathcal{O}_{X_{f}}\left[U\right]\otimes_{\mathcal{O}_{X_{f}}}\mathcal{O}_{X_{f}}\left[T\right],U\mapsto U\otimes1+h^{m-n_{f}}g\otimes T\]
 In other words, the restriction of $\mu_{s}$ to $W_{\gamma}\mid_{X_{f}}$
is simply a translation, twisted by $h^{m-n_{f}}g$. 
\end{rem}
\noindent We also recall the following useful result. 

\begin{prop}
\emph{(}\cite[Proposition 3.12]{DubGDS}\emph{)} \label{TrivialCanonicalClassProp}For
a $GDS$ $q_{\gamma}:W_{\gamma}\rightarrow Z$ defined by a weighted
rooted tree $\gamma=\left(\Gamma,w\right)$, the following are equivalent.

1) $W_{\gamma}$ admits a free $\gaz$-action.

2) All the leaves of $\Gamma$ are at the same level.

3) The canonical sheaf $\omega_{W_{\gamma}}$ of $W_{\gamma}$ is
trivial.
\end{prop}
\begin{example}
In \cite{BML01} Bandman and Makar-Limanov proved that the canonical
sheaf of the surface $V$ of example \ref{Example1} is not trivial
by constructing certain explicit global holomorphic $2$-forms on
$V$. On the other hand, the fine-labelled tree $\mathfrak{g}$ corresponding
to this surface has leaves at levels $1$ and $2$. Therefore, propositions
\ref{Flrt2FwrtTheorem} and \ref{TrivialCanonicalClassProp} implies
the same result. This shows in particular that there is no free $\gaz$-action
on $V$ (see also example \ref{DerOnBMLEx} below). 
\end{example}

\subsection{$\gaz$-actions on the $GDS$'s $q_{\mathfrak{g}}:V_{\mathfrak{g}}\rightarrow Z$ }

\indent

\noindent In this subsection, we prove that every $\gaz$-action
on a $GDS$ $q_{\mathfrak{g}}:V_{\mathfrak{g}}\rightarrow Z$ defined
by a fine-labelled rooted tree extends to a $\gaz$-action on the
ambient space $\ba_{Z}^{d\left(\Gamma\right)+1}=\textrm{Spec}\left(A\left[\Gamma\right]\right)$. 

\begin{example}
\label{DerOnBMLEx} The fine-labelled tree $\mathfrak{g}$ of example
\ref{Example1} corresponds to the Bandman and Makar-Limanov surface
$V\subset\textrm{Spec}\left(A\left[X_{0},X_{e_{0}},X_{e_{1}}\right]\right)$
with equations \[
hX_{e_{0}}=X_{0}\left(X_{0}^{2}-1\right),\quad X_{0}X_{e_{1}}-X_{e_{0}}\left(X_{e_{0}}^{2}-1\right),\quad hX_{e_{1}}=\left(X_{0}^{2}-1\right)\left(X_{e_{0}}^{2}-1\right).\]
 Consider the following locally nilpotent $A$-derivation $\der=\tilde{\der}_{\mathfrak{g},2}$
of $A\left[\Gamma\right]$ given via \[
\left\{ \begin{array}{c}
\der\left(X_{0}\right)=h^{2},\quad\der\left(X_{e_{0}}\right)=h\left(3X_{0}^{2}-1\right)\\
\der\left(X_{e_{1}}\right)=2hX_{0}\left(X_{e_{0}}^{2}-1\right)+2\left(X_{0}^{2}-1\right)\left(3X_{0}^{2}-1\right)X_{e_{0}}\end{array}\right..\]
 In \cite[p. 579]{BML01}, Bandman and Makar-Limanov proved that $\der$
induces a locally nilpotent $A$-derivation $\der_{\mathfrak{g},2}$
of $B_{\mathfrak{g}}$. Note that corresponding $\gaz$-action restricts
to a free action on the complement of the irreducible components \[
C_{f_{1}}=\left\{ h=0,X_{0}=1\right\} \cap V\quad\textrm{and}\quad C_{f_{2}}=\left\{ h=0,X_{0}=-1\right\} \cap V\]
 of $q_{\mathfrak{g}}^{-1}\left(z_{0}\right)$. 
\end{example}
\noindent More generally, the following proposition asserts that
a $GDS$ $q_{\mathfrak{g}}:V_{\mathfrak{g}}\rightarrow Z$ admits
canonical $\gaz$-actions which are the restrictions to $V_{\mathfrak{g}}$
of certain $\gaz$-actions on $\ba_{Z}^{d\left(\Gamma\right)+1}$. 

\begin{prop}
\label{EmbeddedActionProp} For every $m\geq h\left(\Gamma\right)$,
the $A$-derivation $\tilde{\der}_{\mathfrak{g},m}$ of $A_{h}\left[\Gamma\right]$
defined recursively by \[
\tilde{\der}_{\mathfrak{g},m}=h^{m}\frac{\partial}{\partial X_{0}}+h^{-1}\sum_{e\in\mathfrak{P}\left(\Gamma\right)}\tilde{\der}_{\mathfrak{g},m}\left(Q_{e}\left(\mathfrak{g}\right)\right)\frac{\partial}{\partial X_{e}}\]
 restricts to a triangular derivation of $A\left[\Gamma\right]$.
It induces a locally nilpotent $A$-derivation $\der_{\mathfrak{g},m}$
of the $A$-algebra $B_{\mathfrak{g}}=A\left[\Gamma\right]/I_{\mathfrak{g}}$. 
\end{prop}
\begin{proof}
Given a node $e\in\mathbf{P}\left(\Gamma\right)$ at level $k<h\left(\Gamma\right)$,
$Q_{e}\left(\mathfrak{g}\right)$ only involves the variables $X_{0}$
and $X_{e'}$, $e'\in\textrm{Anc}\left(e\right)$ (see definition
\ref{QPolyDefinition}). We conclude recursively that \[
\tilde{\der}_{\mathfrak{g},m}\left(X_{e}\right)=h^{-1}\sum_{e'\in\textrm{Anc}\left(e\right)\cup\left\{ 0\right\} }{\displaystyle \frac{\der Q_{e}\left(\mathfrak{g}\right)}{\der X_{e'}}}\tilde{\der}_{\mathfrak{g},m}\left(X_{e'}\right)\in h^{m-k-1}A{\displaystyle \left[X_{0},\left(X_{e'}\right)_{e'\in\textrm{Anc}\left(e\right)}\right]}.\]
Thus $\tilde{\der}_{\mathfrak{g},m}$ restricts to a triangular $A$-derivation
of $A\left[\Gamma\right]$ as $m\geq h\left(\Gamma\right)$. By construction,
$\tilde{\der}_{\mathfrak{g},m}$ annihilates $\Delta_{0,e}\left(\mathfrak{g}\right)$
for every $e\in\mathfrak{P}\left(\Gamma\right)$. Given $\left(e,e'\right)\in\mathbf{P}\left(\Gamma\right)\times\textrm{Anc}\left(e\right)$,
we deduce from (\ref{SyzigyRelation}) that \[
h\tilde{\der}_{\mathfrak{g},m}\left(\Delta_{e',e}\left(\mathfrak{g}\right)\right)=\tilde{\der}_{\mathfrak{g},m}\left(h\Delta_{e',e}\left(\mathfrak{g}\right)\right)\in I_{\mathfrak{g}}.\]
 Since $B_{\mathfrak{g}}$ is an integral $A$-algebra, we conclude
that $\tilde{\der}_{\mathfrak{g},m}\left(\Delta_{e',e}\left(\mathfrak{g}\right)\right)\in I_{\mathfrak{g}}$
for otherwise $h$ is a zero divisor in $B_{\mathfrak{g}}$. Thus
$\tilde{\der}_{\mathfrak{g},m}\left(I_{\mathfrak{g}}\right)\subset I_{\mathfrak{g}}$
and so, $\tilde{\der}_{\mathfrak{g},m}$ induces a locally nilpotent
$A$-derivation $\der_{\mathfrak{g},m}$on $B_{\mathfrak{g}}$. 
\end{proof}
\noindent In \cite{Fie94}, Fieseler exploited the vector field associated
to a $\bcp$-action to introduce the concept of fixed point order
of a $\bcp$-action along an invariant subscheme. Here we give an
equivalent definition in terms of locally nilpotent derivations. 

\begin{defn}
Let $V=\textrm{Spec}\left(B\right)$ be an integral affine scheme
over a field $k$ of characteristic $0$. Given a nontrivial $\mathbb{G}_{a,k}$-action
$\alpha:\mathbb{G}_{a,k}\times V\rightarrow V$ corresponding to a
locally nilpotent derivation $\der\in\textrm{LND}_{k}\left(B\right)$
and an invariant subscheme $Y\subset V$ with defining ideal $I_{Y}\subset B$,
the \emph{fixed point order $\mu\left(\alpha,Y\right)$ of $\alpha$
along $Y$} is the maximal number $n\geq0$ such that $\der\left(B\right)\subset I_{Y}^{n}B$. 
\end{defn}
\noindent We have the following result.

\begin{lem}
\label{EmbeddedFixedPointLemma} For every $m\geq h\left(\Gamma\right)$
the $\gaz$-action $\alpha_{\mathfrak{g},m}$ on $V_{\mathfrak{g}}$
corresponding to the derivation $\der_{\mathfrak{g},m}$ restricts
to a free action on $V_{\mathfrak{g}}\times_{Z}Z_{*}$. It has fixed
point order $\mu\left(\alpha_{\mathfrak{g},m},C_{e}\right)=m-n_{e}$
along an irreducible component $C_{e}\subset q_{\mathfrak{g}}^{-1}\left(z_{0}\right)$
corresponding to a leaf $e$ of $\Gamma$ at level $n_{e}$. 
\end{lem}
\begin{proof}
Clearly, $\alpha_{\mathfrak{g},m}$ restricts to a free action on
$V_{\mathfrak{g}}\times_{Z}Z_{*}\simeq\textrm{Spec}\left(A_{h}\left[X_{0}\right]\right)$
(see \ref{TrivialBundle}) as $\tilde{\der}_{\mathfrak{g},m}\left(X_{0}\right)=h^{m}$.
With the notation of (\ref{MaximalChainNotation}), the ideal $\tilde{I}_{\mathfrak{g}}\left(e\right)\subset B_{\mathfrak{g}}=A\left[\Gamma\right]/I_{\mathfrak{g}}$
of the curve $C_{e}$ is generated by the image of the ideal \[
I_{\mathfrak{g}}\left(e\right)=\left(I_{\mathfrak{g}},h,X_{0}-lb\left(e_{1}\right),\ldots,X_{e_{n_{e}-2}}-lb\left(e_{n_{e}}\right)\right)\subset A\left[\Gamma\right]\]
 in $B_{\mathfrak{g}}$ (see proposition \ref{ReducedFiberProposition}).
By definition of $\tilde{\der}_{\mathfrak{g},m}$ , \[
\tilde{\der}_{\mathfrak{g},m}\left(h,X_{e_{0}},\ldots,X_{e_{n_{e}-1}}\right)\subset I_{\mathfrak{g}}^{m-n_{e}}\left(e\right)\setminus I_{\mathfrak{g}}^{m-n_{e}+1}\left(e\right).\]
 For $e'\in\mathbf{P}\left(\Gamma\right)\setminus\left(\downarrow e\right)$
there exists an indice $k$, $0\leq k\leq n_{e}-1$ such that $e_{k}$
is the first common ancestor of $e'$ and $e$. Letting $e"=\textrm{Child}\left(e_{k}\right)\cap\left(\downarrow e'\right)$,
we have $Q_{e'}=R_{e''}Q_{e_{k},e'}$, where $Q_{e_{k},e'}$ only
involves the variables corresponding to nodes in $\left[e_{k},\textrm{Par}\left(e'\right)\right]_{\Gamma}$.
Since $\left(X_{e_{k-1}}-lb\left(e_{k+1}\right)\right)$ divides $R_{e''}$
we conclude by induction that \[
\tilde{\der}_{\mathfrak{g},m}\left(X_{e'}\right)=h^{-1}\tilde{\der}_{\mathfrak{g},m}\left(Q_{e'}\right)\in I_{\mathfrak{g}}^{m-k-1}\left(e\right)\subset I_{\mathfrak{g}}^{m-n_{e}}\left(e\right).\]
Thus $\der_{\mathfrak{g},m}\left(B_{\mathfrak{g}}\right)\subset\tilde{I}_{\mathfrak{g}}^{m-n_{e}}\left(e\right)\setminus\tilde{I}_{\mathfrak{g}}^{m-n_{e}-1}\left(e\right)$
and so, $\mu\left(\alpha_{\mathfrak{g},m},C_{e}\right)=m-n_{e}$.
\end{proof}
\noindent This leads to the following description.

\begin{thm}
\label{ActionExtensionProp} Every nontrivial $\gaz$-action on a
$GDS$ $q_{\mathfrak{g}}:V_{\mathfrak{g}}\rightarrow Z$ is induced
by a locally nilpotent derivation $g\der_{\mathfrak{g},m}$ of $B_{\mathfrak{g}}$,
where $m\geq h\left(\Gamma\right)$ and $g\in A\setminus hA$.
\end{thm}
\begin{proof}
Since $A$ is contained in the kernel of $\der_{\mathfrak{g},m}$,
$g\der_{\mathfrak{g},m}$ is again locally nilpotent, whence induces
a nontrivial $\gaz$-action $\alpha$ on $V_{\mathfrak{g}}$. For
every $e\in\textrm{Leaves}\left(\Gamma\right)$, the open subset \[
V_{e}=\left(V_{\mathfrak{g}}\setminus q_{\mathfrak{g}}^{-1}\left(z_{0}\right)\right)\cup C_{e}\]
 is $Z$-isomorphic to $\ba_{Z}^{1}=\textrm{Spec}\left(A\left[U\right]\right)$
(see \ref{GDS2WTree}). By lemma \ref{EmbeddedFixedPointLemma}, $\mu\left(\alpha,C_{e}\right)=m-n_{e}$
, where $n_{e}=l\left(\downarrow e\right)$. We conclude similarly
that $\mu\left(\alpha,q_{\mathfrak{g}}^{-1}\left(z\right)\right)=\textrm{ord}_{z}\left(g\right)$
for every $z\in Z_{*}$. Letting $\gaz=\textrm{Spec}\left(A\left[T\right]\right)$,
this means that the group co-action corresponding to the restriction
of $\alpha$ to $V_{f}$ is given by the $A$-algebras homomorphism
\[
A\left[U\right]\rightarrow A\left[U\right]\otimes_{A}A\left[T\right],U\mapsto U\otimes1+h^{m-n_{e}}g\otimes T.\]
By proposition \ref{Flrt2FwrtTheorem}, the $GDS$ $q_{\gamma}:W_{\gamma}\rightarrow Z$
defined by the weighted rooted tree $\gamma=\left(\Gamma,w\right)$
corresponding to the pair $\left(q_{\mathfrak{g}}:V_{\mathfrak{g}}\rightarrow Z,pr_{X_{0}}:V_{\mathfrak{g}}\rightarrow\ba_{Z}^{1}\right)$
is isomorphic to $V_{\mathfrak{g}}$ via the open embedding $i_{\gamma}:W_{\gamma}\hookrightarrow\textrm{Spec}\left(A\left[\Gamma\right]\right)$.
We deduce from remark \ref{Co-actionRem} that this isomorphism is
equivariant when we equip $V_{\mathfrak{g}}$ with $\alpha$ and $W_{\gamma}$
with the action $\mu_{s}$ corresponding to $h^{m}g\in A$ (see \ref{Function2Action}).
Since every nontrivial $\gaz$-action on $W_{\gamma}$ is of the form
$\mu_{s}$ for some section $s\in H^{0}\left(X_{\gamma},\mathcal{L}_{\gamma}^{\vee}\right)$
corresponding to a regular function $h^{m}g\in A$, where $m\geq h\left(\Gamma\right)$
and $g\in A\setminus hA$, we conclude that every $\gaz$-action on
$V_{\mathfrak{g}}$ is induced by a locally nilpotent derivation $g\der_{\mathfrak{g},m}$
for some $g\in A\setminus hA$. 
\end{proof}
\noindent As a consequence of theorems \ref{Flrt2GDSClassifTheorem}
and \ref{ActionExtensionProp}, we obtain the following result. 

\begin{cor}
Any $GDS$ $q:V\rightarrow Z$ can be embedded into an affine space
$\ba_{Z}^{N}$ in such a way that every $\gaz$-action on $V$ extends
to a $\gaz$-action on $\ba_{Z}^{N}$. 
\end{cor}

\section{$GDS$'s with a trivial Makar-Limanov invariant }

In this section we characterize $GDS$'s $q_{\mathfrak{g}}:V_{\mathfrak{g}}\rightarrow Z$
with a trivial Makar-Limanov invariant. By theorem 7.7 in \cite{DubGDS},
every normal affine surface $S$ with a trivial Makar-Limanov invariant
is a cyclic quotient of a $GDS$ $V$. In case that $S$ is a log
$\mathbb{Q}$-homology planes with a trivial Makar-Limanov invariant,
we construct this $GDS$ $V$ and the cyclic group action explicitly.

\subsection{Embeddings of $GDS$'s with a trivial Makar-Limanov invariant}

\indent

\noindent We recall that an (\emph{oriented}) \emph{comb} of height
$n$ is a rooted tree $\Gamma$ such that for every node $e\in N\left(\Gamma\right)\setminus\textrm{Leaves}\left(\Gamma\right)$
of degree $\deg_{\Gamma}\left(e\right)\geq1$, all but possibly one
of the children of $e$ are leaves of $\Gamma$. This means equivalently
that \begin{equation}
C_{\Gamma}=\left(N\left(\Gamma\right)\setminus\textrm{Leaves}\left(\Gamma\right)\right)=\left\{ e_{0}<\ldots<e_{n-1}\right\} \label{CombNotation}\end{equation}
 is either empty or a chain of length $n-1$. 

\psset{unit=1cm}

\begin{pspicture} (-6.5,-2)(4,0.8)

\psline(-2,0)(3,0)

\psline(-1,0)(-0.8,-1)\psline(-1,0)(-0.6,-0.7)\psline(-1,0)(-0.4,-0.4)

\psline(1,0)(1.2,-1)\psline(1,0)(1.4,-0.7)

\psline(2,0)(2.2,-1)\psline(2,0)(2.4,-0.7)\psline(2,0)(2.6,-0.4)\psline(2,0)(2.8,-0.2)

\rput(-2,0){\textbullet} \rput(-1,0){\textbullet}\rput(0,0){\textbullet}\rput(1,0){\textbullet}

\rput(2,0){\textbullet}\rput(3,0){\textbullet}

\rput(-0.8,-1){\textbullet} \rput(-0.6,-0.7){\textbullet}\rput(-0.4,-0.4){\textbullet}

\rput(1.2,-1){\textbullet} \rput(1.4,-0.7){\textbullet}

\rput(2.8,-0.2){\textbullet} \rput(2.4,-0.7){\textbullet}\rput(2.6,-0.4){\textbullet}\rput(2.2,-1){\textbullet}

\rput(-2,0.2){{\small $e_0$}}

\psframe[linestyle=dotted](-2.2,-0.4)(2.2,0.4)

\rput(0.5,-1.5){\textrm{A comb rooted in } $e_0$.}

\end{pspicture} 

\noindent We have the following result.

\begin{thm}
\label{Comb2GDSThm} A $GDS$ $q_{\mathfrak{g}}:V_{\mathfrak{g}}\rightarrow Z$
has a trivial Makar-Limanov invariant if and only if $\mathfrak{g}=\left(\Gamma,lb\right)$
is a fine-labelled comb.
\end{thm}
\begin{proof}
By proposition \ref{Flrt2FwrtTheorem}, $V_{\mathfrak{g}}$ is $Z$-isomorphic
to the $GDS$ $q_{\gamma}:W_{\gamma}\rightarrow Z$ defined by the
weighted rooted tree $\gamma$ corresponding to the pair $\left(q_{\mathfrak{g}}:V_{\mathfrak{g}}\rightarrow Z,pr_{X_{0}}:V_{\mathfrak{g}}\rightarrow\ba_{Z}^{1}\right)$.
Since $\gamma$ has the same underlying tree $\Gamma$ as $\mathfrak{g}$,
the statement follows from theorem 7.2 in \cite{DubGDS} which asserts
that a $GDS$ $q_{\gamma}:W_{\gamma}\rightarrow Z$ has a trivial
Makar-Limanov invariant if and only if $\Gamma$ is a comb.
\end{proof}
\begin{enavant} \label{NewEquations} Given a nontrivial fine-labelled
comb $\mathfrak{g}=\left(\Gamma,lb\right)$ of height $n$, the ideal
$I_{\mathfrak{g}}\subset A\left[\Gamma\right]=A\left[X_{0}\right]\left[X_{e_{0}},\ldots,X_{e_{n-1}}\right]$
of $V_{\mathfrak{g}}$ is generated by the polynomials\begin{equation}
\left\{ \begin{array}{lcll}
\Delta_{0,e_{i}}\left(\mathfrak{g}\right) & = & hX_{e_{i}}-Q_{e_{i}}\left(\mathfrak{g}\right) & 0\leq i\leq n-1\\
\Delta_{e_{j},e_{i}}\left(\mathfrak{g}\right) & = & \left(X_{e_{j-1}}-lb\left(e_{j+1}\right)\right)X_{e_{i}}-X_{e_{j}}Q_{e_{j,}e_{i}}\left(\mathfrak{g}\right) & 0\leq j<i\leq n-1,\end{array}\right.\label{MLEquationCase}\end{equation}
 \noindent These equations can be put in the following more convenient
form. Letting \[
P_{i}\left(T\right)=\left(T-\lambda_{i,1}\right)\tilde{P}_{i}\left(T\right)=\left(T-lb\left(e_{i}\right)\right)S_{e_{i-1}}^{\left\{ e_{i}\right\} }\left(\mathfrak{g}\right)\left(T\right)\in\bc\left[T\right],\quad1\leq i\leq n,\]
 it is easily seen that $V_{\mathfrak{g}}$ is isomorphic to the surface
\[
V_{P_{1},\ldots,P_{n}}\subset\ba_{\bc}^{n+2}=\textrm{Spec}\left(\bc\left[x\right]\left[y_{1},\ldots,y_{n+1}\right]\right)\]
 with equations \[
\left\{ \begin{array}{lcll}
xy_{i+1} & = & \left({\displaystyle \prod_{k=1}^{i-1}}\tilde{P}_{k}\left(y_{k}\right)\right)P_{i}\left(y_{i}\right) & 1\leq i\leq n\\
\left(y_{j-1}-\lambda_{j-1,1}\right)y_{i+1} & = & y_{j}\left({\displaystyle \prod_{k=j}^{i-1}}\tilde{P}_{k}\left(y_{k}\right)\right)P_{i}\left(y_{i}\right) & 2\leq j\leq i\leq n\end{array}\right.,\]
 where, by convention, ${\displaystyle {\displaystyle \prod_{k=j}^{i-1}}}\tilde{P}_{k}\left(y_{k}\right)=1$
if $i=j$ . In particular, if $\mathfrak{g}$ is a fine-labelled comb
of height $1$ then $V_{\mathfrak{g}}$ is the ordinary Danielewski
surface with equation $xy_{2}-P_{1}\left(y_{1}\right)=0$. 

\end{enavant}

\begin{rem}
Note that for every collection of nonconstant polynomials $P_{1},\ldots,P_{n}$
with simple roots, the surface $V=V_{P_{1},\ldots,P_{n}}$ admits
two natural $\ba^{1}$-fibrations $pr_{x}:V\rightarrow\textrm{Spec}\left(\bc\left[x\right]\right)$
and $pr_{y_{n+1}}:V\rightarrow\textrm{Spec}\left(\bc\left[y_{n+1}\right]\right)$.
By construction $pr_{x}:V\rightarrow\textrm{Spec}\left(\bc\left[x\right]\right)$
is a $GDS$ over $Z$. However, it may happen that $pr_{y_{n+1}}:V\rightarrow\textrm{Spec}\left(\bc\left[y_{n+1}\right]\right)$
is a not a $GDS$, due to the fact that the fiber $pr_{y_{n+1}}^{-1}\left(0\right)$
is not reduced (see example \ref{NonConjEx} below) 
\end{rem}
\noindent Proposition \ref{TrivialCanonicalClassProp} implies the
following result due to Bandman and Makar-Limanov \cite[Theorem 3]{BML01}.

\begin{thm}
\label{MLGDSThm} For a $GDS$ $q:V\rightarrow Z$ with a trivial
Makar-Limanov invariant, the following are equivalent.

1) $V$ admits a free $\gaz$-action.

2) The canonical sheaf $\omega_{V}$ of $V$ is trivial.

3) $V$ is isomorphic to an ordinary Danielewski surface $V_{P,1}\subset\ba_{\bc}^{3}=\textrm{Spec}\left(\bc\left[x,y,z\right]\right)$
with the equation $xz-P\left(y\right)=0$ for a certain nonconstant
polynomial $P$ with simple roots.
\end{thm}
\begin{proof}
By virtue of proposition \ref{TrivialCanonicalClassProp}, (1) and
(2) are equivalent. Clearly, (3)$\Rightarrow$(2). So it remains to
show that (2)$\Rightarrow$(3). By theorem \ref{Flrt2GDSClassifTheorem},
there exists a fine-labelled rooted tree $\mathfrak{g}=\left(\Gamma,lb\right)$
and a $Z$-isomorphism $\psi:V\stackrel{\sim}{\rightarrow}V_{\mathfrak{g}}$.
By theorem \ref{MLGDSThm}, $\mathfrak{g}$ is a comb as $V$ has
a trivial Makar-Limanov invariant. Since $\omega_{V}$ is trivial,
proposition \ref{TrivialCanonicalClassProp} implies that all the
leaves of $\Gamma$ are at the same level $n\geq0$. Therefore $\Gamma$
has the following structure:

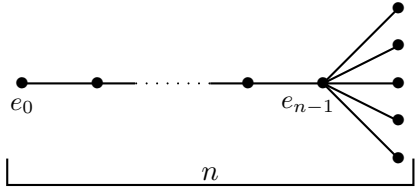
\psset{unit=1cm}

\begin{pspicture}(-4.5,-2)(6,1.5)

\psline(0,0)(1.5,0)

\psline[linestyle=dotted](1.5,0)(2.5,0)

\psline(2.5,0)(4,0)

\psline(4,0)(5,1)\psline(4,0)(5,0.5)\psline(4,0)(5,0)\psline(4,0)(5,-0.5)\psline(4,0)(5,-1)

\rput(0,0){\textbullet}\rput(1,0){\textbullet}\rput(3,0){\textbullet}

\rput(4,0){\textbullet}

\rput(5,0){\textbullet}\rput(5,0.5){\textbullet}\rput(5,1){\textbullet}\rput(5,-0.5){\textbullet}

\rput(5,-1){\textbullet}

\rput(0,-0.3){{\small $e_0$}}

\rput(3.8,-0.3){{\small $e_{n-1}$}}

\pnode(-0.2,-1 ){C}\pnode(5.2,-1){D}\ncbar[angle=-90]{C}{D}

\rput(2.5,-1.2){$n$}

\end{pspicture}

\noindent If $n=0$ then $V_{\mathfrak{g}}\simeq\ba_{Z}^{1}=\textrm{Spec}\left(A\left[X_{0}\right]\right)$.
Otherwise, letting $X_{e_{-1}}=X_{0}$, equations (\ref{MLEquationCase})
simplify as \[
\left\{ \begin{array}{lcll}
hX_{e_{i}} & = & X_{e_{i-1}}-lb\left(e_{i+1}\right) & 0\leq i\leq n-2\\
hX_{e_{n-1}} & = & \prod_{f\in\textrm{Child}\left(e_{n-1}\right)}\left(X_{e_{n-2}}-lb\left(f\right)\right)\end{array}\right.,\]
 and so, the statement follows.
\end{proof}
\noindent By \cite{ML90}, any two free $\bcp$-actions on an ordinary
Danielewski surface $V\subset\textrm{Spec}\left(\bc\left(x,y,z\right)\right)$
with equation $xz-P\left(y\right)=0$ are conjugated under the action
of the automorphism group $\textrm{Aut}\left(V\right)$ of $V$. In
particular, every two $\ba^{1}$-fibrations $q:V\rightarrow Z$ and
$\tilde{q}:V\rightarrow\tilde{Z}$ are conjugated, in the sense that
there exists an isomorphism $\bar{\theta}:Z\stackrel{\sim}{\rightarrow}\tilde{Z}$
and an automorphism $\theta:V\stackrel{\sim}{\rightarrow}V$ of $V$
such that $\bar{\theta}\circ q=\tilde{q}\circ\theta$. In \cite{DaRus02},
Daigle and Russel give examples of nonconjugated $\ba^{1}$-fibrations
on a log $\mathbb{Q}$-homology plane with a trivial Makar-Limanov
invariant. In case that $q:V\rightarrow Z$ is a $GDS$, the following
example shows that there may also exist nonconjugated $\ba^{1}$-fibrations
on $V$. 

\begin{example}
\label{NonConjEx} By example \ref{Example1}, the Bandman and Makar-Limanov
surface $V\subset\textrm{Spec}\left(\bc\left[x,y,z,u\right]\right)$
with equations\[
xz=y\left(y^{2}-1\right),\quad yu=z\left(z^{2}-1\right),\quad xu=\left(y^{2}-1\right)\left(z^{2}-1\right)\]
 is a $GDS$ via $pr_{x}:V\rightarrow Z=\textrm{Spec}\left(\bc\left[x\right]\right)$.
It is also a $GDS$ via $pr_{u}:V\rightarrow\tilde{Z}=\textrm{Spec}\left(\bc\left[u\right]\right)$.
The $\bc$-algebra involution \[
\bc\left[x,u\right]\left[y,z\right]\rightarrow\bc\left[x,u\right]\left[y,z\right],x\leftrightarrow u,y\leftrightarrow z,\]
 induces an automorphism $\theta:V\stackrel{\sim}{\rightarrow}V$.
Letting $\bar{\theta}:Z\stackrel{\sim}{\rightarrow}\tilde{Z}$ be
the isomorphism corresponding to the induced isomorphism $\bc\left[u\right]\stackrel{\sim}{\rightarrow}\bc\left[x\right]$,
we conclude that $\bar{\theta}\circ pr_{x}=pr_{u}\circ\theta$.

Similarly, the surface $V'\subset\bc\left[x,y,z,u\right]$ with equations
\[
xz=y\left(y-1\right),\quad yu=z^{2},\quad xu=\left(y-1\right)z\]
 is a $GDS$ $pr_{x}:V'\rightarrow Z=\textrm{Spec}\left(\bc\left[x\right]\right)$.
The natural second $\ba^{1}$-fibration $pr_{u}:V'\rightarrow\tilde{Z}$
restricts to the trivial line bundle over $\tilde{Z}\setminus\left\{ 0\right\} $
but the fiber\[
pr_{u}^{-1}\left(0\right)\simeq\textrm{Spec}\left(\bc\left[x,y,z\right]/\left(z^{2},\left(y-1\right)z,xz-y\left(y-1\right)\right)\right)\]
 is not reduced. Therefore, there is no automorphism $\theta:V'\stackrel{\sim}{\rightarrow}V'$
satisfying $\bar{\theta}\circ pr_{x}=pr_{u}\circ\theta$. 
\end{example}

\subsection{Log $\mathbb{Q}$-homology planes with a trivial Makar-Limanov invariant}

\indent

\noindent A log $\mathbb{Q}$-homology plane is a normal affine surface
$S$ with log-terminal singularities and with vanishing singular homology
groups $H_{i}\left(S,\mathbb{Q}\right)$ for all $i>0$. In this subsection
we prove the following result (compare with \cite[Theorem 3.4]{MiMa03}).

\begin{thm}
\label{LogQHomThm} Every log $\mathbb{Q}$-homology plane $S\not\simeq\ba_{\bc}^{2}$
with a trivial Makar-Limanov invariant is isomorphic to the quotient
of an ordinary Danielewski surface $V\subset\textrm{Spec}\left(\bc\left[x,t,z\right]\right)$
with equation $xz=t^{n}-1$ by a $\mathbb{Z}_{m}$-action $\left(x,t,z\right)\mapsto\left(\varepsilon x,\varepsilon^{q}t,\varepsilon^{-1}z\right)$,
where $\varepsilon$ is a primitive $m$-th root of unity, $n$ divides
$m$ and $\gcd\left(q,m/n\right)=1$. 
\end{thm}
\noindent The proof is given in \ref{QhomStep1}-\ref{QhomEnd} below.
Since $S$ has a trivial Makar-Limanov invariant, it admits an $\ba^{1}$-fibration
$\rho:S\rightarrow Y\simeq\textrm{Spec}\left(\bc\left[y\right]\right)$
such that all but possibly one closed fiber, say $\rho^{-1}\left(y_{0}\right)$,
where $y_{0}=\left(y\right)\in Y$, are isomorphic to $\ba_{\bc}^{1}$
(see e.g. \cite[Proposition 2.15]{Dub02}). In turn, lemma 1.2 in
\cite{Fie94} and \cite[2.5]{DubGDS} imply that $\rho$ restricts
to a trivial line bundle over $Y_{*}=Y\setminus\left\{ y_{0}\right\} $.
By Theorem 4.3.1 in \cite{MiyBook} p.231, $\rho^{-1}\left(y_{0}\right)=mC$,
where $m\geq1$ and $C\simeq\ba_{\bc}^{1}$. In particular, there
exists $k$, $0\leq k\leq m-1$, such that the canonical sheaf $\omega_{S}$
of $S$ is isomorphic to $\mathcal{O}_{S}\left(kC\right)$. 

\begin{enavant} \label{QhomStep1} If $m=1$ then $S\simeq\ba_{Y}^{1}$
(see e.g. \cite[Theorem 4.3.1 p.231]{MiyBook}). Otherwise, if $m\geq2$
then we let $\theta:Z\simeq\ba_{\bc}^{1}\rightarrow Y$ be a Galois
covering of order $m$, unramified over $Y_{*}$ and totally ramified
over $y_{0}$. Up to an isomorphism, we can suppose that $\theta$
is induced by the $\bc$-algebras homomorphism $\theta^{*}:\bc\left[y\right]\rightarrow\bc\left[x,y\right]/\left(x^{m}-y\right)\simeq\bc\left[x\right]$,
where the Galois group $\mathbb{Z}_{m}$ of $m$-th roots of unity
acts on $\bc\left[x,y\right]$ via $\varepsilon\cdot f\left(x,y\right)=f\left(\varepsilon^{-1}x,y\right)$.
The group $\mathbb{Z}_{m}$ acts on the normalization $V$ of the
fiber product $S\times_{Y}Z$ and $S\simeq V/\mathbb{Z}_{m}$. This
surface $V$ inherits an $\ba^{1}$-fibration $q:V\rightarrow Z$
with reduced fibers, induced by the projection $pr_{Z}:S\times_{Y}Z\rightarrow Z$.
Letting $z_{0}=\left(x\right)\in Z$, the group $\mathbb{Z}_{m}$
acts transitively on the set of irreducible components $C_{1},\ldots,C_{n}$
of $q^{-1}\left(z_{0}\right)$, where $n$ divides $m$ (see e.g.
\cite[Theorem 1.7]{Fie94}). Since $\rho:S\rightarrow Y$ restricts
to a trivial line bundle over $Y_{*}$, we conclude that $q:V\rightarrow Z$
restricts to a trivial line bundle over $Z_{*}=Z\setminus\left\{ z_{0}\right\} $
and so, $q:V\rightarrow Z$ is a $GDS$ over $Z$. By construction,
the quotient morphism $\phi:V\rightarrow S$ is étale over $S_{reg}$,
whence is étale in codimension $1$ as $S$ is normal. Thus the canonical
sheaf $\omega_{V}$ of $V$ is isomorphic to $\phi^{*}\left(\omega_{S}\right)\simeq q^{*}\left(\mathcal{O}_{Z}\left(kz_{0}\right)\right)$,
hence is trivial. Moreover, we deduce from \cite{Vas69} that the
Makar-Limanov invariant of $V$ is trivial.

\end{enavant}

\begin{enavant} \label{QhomEnd} By theorem \ref{MLGDSThm}, $V$
is isomorphic to an ordinary Danielewski surface $V_{P,1}\subset\textrm{Spec}\left(\bc\left[x,t,z\right]\right)$
with the equation $xz-P\left(t\right)=0$ for a certain nonconstant
polynomial $P$ with $n$ simple roots $t_{1},\ldots,t_{n}\in\bc$.
Since $\mathbb{Z}_{m}$ acts transitively on the set of irreducible
components of $q^{-1}\left(z_{0}\right)$ and $t$ is locally constant
on $q^{-1}\left(z_{0}\right)$ with the value $t_{i}$ on $C_{i}$,
$1\leq i\leq n$, we can suppose that $t_{i+1}=\varepsilon^{q}t_{i}$
where $\varepsilon$ is a primitive $m$-th root of unity and $\gcd\left(q,m/n\right)=1$.
Up to an isomorphism, we can further suppose that $t_{1}=1$ (see
e.g. \cite{ML90}). The canonical morphism $pr_{x,t}:V\rightarrow\ba_{\bc}^{2}=\textrm{Spec}\left(\bc\left[x,t\right]\right)$
expresses $V$ as the affine modification \cite{KaZa99} of $\ba_{\bc}^{2}$
with the locus $\left(I,f\right)=\left(\left(x,P\left(t\right)\right),x\right)$.
Clearly, $pr_{x,t}$ is $\mathbb{Z}_{m}$-equivariant when we equip
$\ba_{\bc}^{2}$ with the $\mathbb{Z}_{m}$-action $\left(x,t\right)\mapsto\left(\varepsilon x,\varepsilon^{q}t\right)$.
As $\Gamma\left(V,\mathcal{O}_{V}\right)\simeq\bc\left[x,t,P\left(t\right)/x\right]$,
we conclude that $V$ is isomorphic to the ordinary Danielewski surface
with equation $xz=t^{n}-1$, equipped with the $\mathbb{Z}_{m}$-action
$\left(x,t,z\right)\mapsto\left(\varepsilon x,\varepsilon^{q}t,\varepsilon^{-1}z\right)$. 

\end{enavant}

\begin{rem}
Note that $S\simeq V/\mathbb{Z}_{m}$ is nonsingular if and only if
$n=m$. On the other hand, if $n=1$ then $S$ is isomorphic to a
cyclic quotient of $\ba^{2}=\textrm{Spec}\left(\bc\left[x,t\right]\right)$
by the $\mathbb{Z}_{m}$-action $\left(x,t\right)\mapsto\left(\varepsilon x,\varepsilon^{q}t\right)$,
where $\gcd\left(q,m\right)=1$. In particular, every affine toric
surface nonisomorphic to $\bc\times\bc^{*}$ has a trivial Makar-Limanov
invariant, see also \cite[Proposition 2.11]{Dub02} for a geometrical
proof.
\end{rem}
\begin{cor}
Every log $\mathbb{Q}$-homology plane with a trivial Makar-Limanov
invariant admits a nontrivial $\bc^{*}$-action.
\end{cor}
\begin{proof}
Indeed, every ordinary Danielewski surface $V_{P,1}\subset\ba_{\bc}^{3}=\textrm{Spec}\left(\bc\left[x,t,z\right]\right)$
with equation $xz-P\left(t\right)=0$ admits a nontrivial $\bc^{*}$-action
induced by the $\bc^{*}$-action on $\ba_{\bc}^{3}$ $\lambda\cdot\left(x,t,z\right)=\left(\lambda x,y,\lambda^{-1}z\right)$,
$\lambda\in\bc^{*}$. In case that $P\left(t\right)=t^{n}-1$, this
action commutes with the $\mathbb{Z}_{m}$-action $\left(x,t,z\right)\mapsto\left(\varepsilon x,\varepsilon^{q}t,\varepsilon^{-1}z\right)$
and so, the quotient $S=V/\mathbb{Z}_{m}$ inherits a nontrivial $\bc^{*}$-action. 
\end{proof}
\bibliographystyle{amsplain}

\end{document}